\documentclass[12pt]{amsart}
\usepackage{amssymb, amsthm, enumerate, enumitem, graphicx, amscd, amsmath, amssymb, amsfonts, mathrsfs, xcolor, tikz, tikz-cd, tkz-euclide, comment}
\usepackage[mathcal]{euscript}

\numberwithin{equation}{section}
\newtheorem{theorem}{Theorem}
\numberwithin{theorem}{section}
\newtheorem{definition}[theorem]{Definition}
\newtheorem{notation}[theorem]{Notation}
\newtheorem{lemma}[theorem]{Lemma}

\title[Rhaly matrices]{Boundedness, compactness and Schatten class for Rhaly matrices}

\author[C. Bellavita]{Carlo Bellavita}
\email{carlo.bellavita@gmail.com}
\email{carlobellavita@ub.edu}
\address{Departament of Matem\'atica i Inform\'atica, Universitat de Barcelona, Gran Via 585, 08007 Barcelona, Spain.}

\author[E. Dellepiane]{Eugenio Dellepiane}
\email{dellepianeeugenio@gmail.com}
\email{eugenio.dellepiane@polito.it}
\address{Dipartimento di Scienze Matematiche “Giuseppe Luigi Lagrange”, Politecnico di Torino, Corso Duca degli Abruzzi 24, 10129 Torino, Italy}

\author[G. Stylogiannis]{Georgios Stylogiannis}
\email{g.stylog@gmail.com}
\email{stylog@math.auth.gr}
\address{Department of Mathematics, Aristotle University of Thessaloniki, 54124, Greece.}

\thanks{The first two authors are members of Gruppo Nazionale per l’Analisi Matematica, la Probabilit\`a e le loro Applicazioni (GNAMPA) of Istituto Nazionale di Alta Matematica (INdAM).
The first author was partially supported by PID2021-123405NB-I00 by
the Ministerio de Ciencia e Innovaci\'on. 
The third author was partially supported by the Hellenic Foundation for Research and Innovation (H.F.R.I.) under the '2nd Call for H.F.R.I. Research Projects to support Faculty Members \& Researchers' (Project Number: 4662).}

\subjclass{47B10,
30H25, 30H99}
\keywords{Rhaly matrices; Schatten classes; Besov spaces; Hadamard multipliers; Weighted harmonically Dirichlet spaces.}

\begin{document}
\begin{abstract}
In this article we present new proofs for the boundedness and the compactness on $\ell^2$ of the Rhaly matrices, also known as terraced matrices. We completely characterize when such matrices belong to the Schatten class $\mathcal{S}^q(\ell^2)$, for $1<q<\infty$. Finally, we apply our results to study the Hadamard multipliers in weighted Dirichlet spaces, answering a question left open by Mashreghi--Ransford.
\end{abstract}
\maketitle
\tableofcontents 

\section*{Introduction}
Given $(\alpha_n)_{n\in\mathbb{N}}$ a sequence of complex numbers, we introduce the Rhaly operator $R_\alpha$ on $\ell^2$, the Hilbert space of square-summable functions $f\colon\mathbb{N}\to\mathbb{C}$,
in the following way:
\[\big(R_\alpha f\big)(k)=\alpha_k\sum_{j=0}^k f(j).\]
During this work, we will also denote the elements of $\ell^2$ as square-summable infinitely long vectors $x=(x_0,x_1,x_2\ldots)$. Using the standard orthonormal basis of $\ell^2$, that is, $\{\textbf{e}^j\}_{j\in\mathbb{N}}$, where $\textbf{e}^j(k) = \delta_{k,j}$ is the Kronecker delta, we represent the action of $R_\alpha$ with the infinite matrix
\begin{equation*}\label{E:terracedmatrix}
\begin{pmatrix}
\alpha_0 & 0 & 0 & 0 & 0 & \cdots\\
\alpha_1 & \alpha_1 & 0 & 0 & 0 & \cdots\\
\alpha_2 & \alpha_2 & \alpha_2 & 0 & 0 & \cdots \\
\alpha_3 & \alpha_3 & \alpha_3 & \alpha_3 & 0 & \cdots\\
\vdots & \vdots & \vdots & \vdots &\vdots & \ddots
\end{pmatrix} .
\end{equation*}
Rhaly matrices arise as the natural generalization of the classical Ces\`{a}ro operator. Choosing 
\[\alpha_k=\frac{1}{k+1},\qquad k\in\mathbb{N},\]
then $R_\alpha$ coincides with the classical Ces\`{a}ro operator $\mathcal{C}$, that assigns to each $f\in\ell^2$ the averages
\[\mathcal{C}f(k)=\frac{1}{k+1}\sum_{j=0}^kf(j), \qquad k\in\mathbb{N}.\]
In \cite{10.1112/blms/21.4.399}, Rhaly showed that whenever the limit $L = \lim_{n \to \infty}(n + 1)|\alpha_n|$ exists, it provides a practical test for the boundedness of $R_\alpha$. Subsequently, Rhaly fixed his attention on the spectral properties, trying to understand which characteristics of the Ces\`{a}ro operator hold also for more general sequences $(\alpha_n)_{n\in\mathbb{N}}$. In \cite{10.1112/blms/21.4.399} he examined the conditions under which $R_\alpha$ is hyponormal and, among the other works, Rhaly research extended to the study of posinormality, see \cite{2} and \cite{10.4134/BKMS.2009.46.1.117}. More recently, Gallardo-Guti\'errez and Partington \cite{GallardoGutierrez2024RhalyOM} have deepened the study of the spectrum of $R_\alpha$, already started by Rhaly in \cite{10.1112/blms/21.4.399}.

Since $R_\alpha$ is an infinite matrix, the first characterizations of its boundedness and compactness were achieved through linear algebraic methods, as detailed in Theorems 9.2 and 9.4 of \cite{grosse1998blocking}.

 Later, the spaces of analytic functions appeared in the study of Rhaly matrices. For example, in \cite{Galanopoulos2022} and \cite{Galanopoulos2023}, Galanopoulos, Girela and Merchán considered the case when $(\alpha_k)_k$ is the sequence of the moments of a finite complex Borel measure $\mu$ on $[0,1]$, that is, 
\[\alpha_k=\mu_k=\int_{[0,1]} t^k\operatorname{d}\!\mu(t).\]
They showed that, with the usual identification between the Hardy space $H^2$ and $\ell^2$, if $g(z)=\sum_n b_n z^n$ with $b=(b_0,b_1,\ldots)\in\ell^2$, then the integral operator on $H^2$
\[\mathcal{C}_\mu g(z)=\int_{[0,1]}\frac{g(tz)}{1-tz}\,\operatorname{d}\!\mu(t)\]
can be written as the power series 
\[\mathcal{C}_\mu g(z)=\sum_{k \geq 0} (R_\alpha b)(k) \, z^k,\]
and in particular it holds the identity $\|\mathcal{C}_\mu\|_{H^2\to H^2}=\|R_\alpha\|_{\ell^2 \to \ell^2}$. 
We also cite \cite{Galanopoulos2024}, where the action of Rhaly matrices on non-Hilbert spaces of sequences is studied, and \cite{Lin2018RhalyOO}, where the authors considered Rhaly matrices on  certain  Fock spaces.

In this article, we give original characterizations and new proofs for the boundedness and the compactness of $R_\alpha$. Also, we completely describe the $q$-Schatten class membership for Rhaly matrices, for $1<q<\infty$, and we characterize this property in terms of spaces of analytic functions. This is the first complete result concerning Schatten class properties for Rhaly matrices. Finally, we show that Rhaly matrices can be used to study the Hadamard multipliers on superharmonically weighted Dirichlet spaces. 


\section{Main results}\label{S:mainresults}
In order to analyze the properties of $R_\alpha$, we introduce some spaces of analytic functions on the unit disk.

\subsection{Spaces of analytic functions and notation}
Given a parameter $1<q<\infty$, we define the Besov space $\mathcal{B}_{\frac{1}{2}}^{2,q}$ as
\[\mathcal{B}_{\frac{1}{2}}^{2,q}=\{f \in \operatorname{Hol}(\mathbb{D}) \colon \int_0^1(1-r)^{\frac{q}{2}-1}\biggl(\int_{0}^{2\pi}|f'(re^{i\theta})|^2 \,\operatorname{d}\!\theta\biggr)^{\frac{q}{2}}r\operatorname{d}\!r<\infty\},\]
where $\operatorname{Hol}(\mathbb{D})$ denotes the space of holomorphic functions on the unit disk $\mathbb{D}$. 
Notice that, for $q=2$, we recover the classical Dirichlet space~$\mathcal{D}$
\[
\mathcal{B}_{\frac{1}{2}}^{2,2} = \mathcal{D}=\{f \in \operatorname{Hol}(\mathbb{D}) \colon \int_{\mathbb{D}} |f'(z)|^2 \,\operatorname{d}\!A(z)<\infty\},
\]
where $\operatorname{d}\!A$ is the two-dimensional Lebesgue measure normalized on $\mathbb{D}$. 
For general properties of this class of spaces, see \cite[Chapter 5]{pavlovi2014function}. In particular, we use the characterization \cite[Corollary 5.2]{pavlovi2014function} obtained in terms of the \emph{dyadic difference operator} $\Delta_n$. Let $Z_{n}=\{k\in\mathbb{N}: 2^{n}\leq k< 2^{n+1}\}$ for $n\geq1$ and $Z_{0}=\{0,1\}$. Given an analytic function $f(z)=\sum_{n\geq 0}a_n z^n$, we define 
\[    \Delta_{n}f(z)=\sum_{j \in Z_{n}} a_jz^j.\]
According to Theorem 5.1 in \cite{pavlovi2014function}, for $1<q<\infty$, we have that $f\in\mathcal{B}_{\frac{1}{2}}^{2,q}$ if and only if 
\begin{equation*}
    \sum_{n\geq 0}\bigl( 2^n \|\Delta_nf\|_{H^2}^2\bigr)^{\frac{q}{2}}<\infty.
\end{equation*}

For $q=\infty$, we define
\[    \mathcal{B}_{\frac{1}{2}}^{2,\infty}=\{ f \in \operatorname{Hol}(\mathbb{D}) \colon \sup_{0\leq r<1}(1-r)^{\frac{1}{2}}\biggl(\int_{0}^{2\pi}|f'(re^{i\theta})|^2 \,\operatorname{d}\!\theta\biggr)^{\frac{1}{2}}<\infty\}.\]
This extremal Besov space is equivalent to the \emph{mean Lipschitz} space
\[\Lambda^2_\frac{1}{2}=\{f\in L^2(\mathbb{T})\colon \sup_{t\in[-\pi,\pi]} |t|^{-\frac{1}{2}} \|f_t-f\|_{L^2} <\infty\},\]
where $f_t(e^{i\theta})=f(e^{i(\theta-t)}).$ In particular, it follows from \cite[Proposition 1.3]{Bourdon_Shapiro_Sledd_1989} that a function $f\in H^2(\mathbb{D)}$ belongs to $\mathcal{B}_{\frac{1}{2}}^{2,\infty}$ if and only if its boundary value function belongs to $\Lambda^2_{\frac{1}{2}}$. We also define the \emph{small} mean Lipschitz space
\[   \lambda^2_\frac{1}{2}=\{f\in \Lambda^{2}_{\frac{1}{2}} \colon \lim_{t\to 0} |t|^{-\frac{1}{2}} \|f_t-f\|_{L^2} =0\}.\]
As usual, we identify a function with its boundary values.  In \cite[Corollary 3.2]{Bourdon_Shapiro_Sledd_1989}, it is shown  that an analytic function $f=\sum_{n\geq 0} a_n z^n$ belongs to $\Lambda^2_{\frac{1}{2}}$ if and only if 
\[\sup_{n} 2^n \|\Delta_n f\|_{H^2}^2 <\infty,\]
and that it belongs to $\lambda^2_{\frac{1}{2}}$ if and only if 
\[\lim_{n} 2^n \|\Delta_n f\|_{H^2}^2 =0.\] See  \cite{Bourdon_Shapiro_Sledd_1989} for properties of mean Lipschitz spaces.

We also recall some well-known definitions and properties, needed for the notation of this article. We denote $\mathcal{B}(\ell^2)$ the space of bounded linear operators on $\ell^2$ and $\|\cdot\|$ the operator norm, that is,
\[\|T\|=\sup_{\|f\|_2= 1 }\|Tf\|_2,\]
where $\|\cdot\|_2$ is the standard $\ell^2$-norm.  We denote by $\mathcal{F}$ the set of finite rank operators. The \emph{essential norm} of $T$ is defined as
\[\|T\|_e=\inf_{P\in \mathcal{F}}\|T-P\|.\]
It is clear that $\|\cdot\|_e\leq \|\cdot\|$.
If $\mathcal{F}_n$ is the set of all linear operators on $\ell^2$  whose range has dimension less than or equal to $n$, the  approximation numbers are defined as
\[a_{n+1}(T)=\inf_{P\in \mathcal{F}_n}\|T-P\|.\]    
Note that $a_1(T)=\|T\|$, the sequence $(a_n(T))_n$ is decreasing and
\[\|T\|_e=\lim_{n \to \infty}a_n(T).\]
An operator $K$ is compact if and only if $\lim_{n \to \infty} a_n(K)=0$. Also, given any compact operator $K$, then for every $T\in\mathcal{B}(\ell^2)$ one has that $\|T\|_e=\|T+K\|_e$. Finally, given $1<q<\infty$, we say that the operator $T$ belongs to the Schatten class $\mathcal{S}^q(\ell^2)$ if and only if $(a_n(T))_{n\in\mathbb{N}} \in \ell^q$. For more details, we refer to \cite[Chapter 1]{zhu2007operator}.

\subsection{Main theorems on Rhaly matrices}
In this article, we prove the following theorems establishing connections between the Rhaly matrix $R_\alpha$ and the Besov spaces. Notice that, if the Rhaly matrix $R_\alpha$ defines a linear operator $\ell^2\to\ell^2$, then necessarily $\alpha=R_\alpha \textbf{e}^0 \in\ell^2$. Given $\alpha\in\ell^2$, we denote $f_\alpha=\sum_{n=0}^\infty \alpha_n z^n \in H^2$ the associated analytic function. 

\begin{theorem}\label{T:main1} The following statements are equivalent.
\begin{enumerate}[label=(\roman*)]
    \item $R_\alpha\colon\ell^2\to\ell^2$ is bounded.
    \item $f_\alpha$ belongs to $\Lambda^2_{\frac{1}{2}}$.
    \item It holds that  \[\sup_n 2^n \sum_{j=2^n}^{2^{n+1}-1}|\alpha_j|^2 <\infty.\]
\end{enumerate}
\end{theorem}

\begin{theorem}\label{T:main2} The following  statements  are equivalent.
\begin{enumerate}[label=(\roman*)]
    \item $R_\alpha\colon\ell^2\to\ell^2$ is compact.
    \item $f_\alpha$ belongs to $\lambda^2_{\frac{1}{2}}$.
    \item  It holds that \[\lim_{n\to \infty} 2^n \sum_{j=2^n}^{2^{n+1}-1}|\alpha_j|^2 =0.\]
\end{enumerate}
\end{theorem}

We point out that, to our knowledge, $q$-Schatten class properties for Rhaly matrices when $q\neq 2$ have never appeared in the literature before. 

\begin{theorem}\label{T:main3} 
Let $1<q<\infty$. The following statements are equivalent.
\begin{enumerate}[label=(\roman*)]
    \item \label{T:RaSchatteni} $R_\alpha$ belongs to the Schatten class $\mathcal{S}^q(\ell^2)$.
    \item \label{T:RaSchattenii} $f_\alpha$ belongs to $\mathcal{B}^{2,q}_{\frac{1}{2}}$.
    \item \label{T:RaSchatteniii} It holds that \[\sum_{n=0}^\infty \bigg( 2^n \sum_{j=2^n}^{2^{n+1}-1}|\alpha_j|^2 \bigg)^\frac{q}{2}<\infty.\]
\end{enumerate}
\end{theorem}


\subsection{Hadamard multipliers in weighted Dirichlet spaces}
In Section \ref{S:section6}, we consider another class of matrices that has been recently introduced by Mashreghi and Ransford in \cite{MR4032210}, in connection with the Hadamard multipliers on superharmonically weighted Dirichlet spaces. Given two formal power series $f(z)=\sum_n a_n z^n, g(z)=\sum_n b_n z^n$, their Hadamard product is defined as
\[f\ast g(z) = \sum_n a_n b_n z^n.\]
If $\omega$ is a positive superharmonic function on $\mathbb{D}$, then the associated superharmonically weighted Dirichlet space is 
\[\mathcal{D}_\omega=\{f\in\operatorname{Hol}(\mathbb{D})\colon \int_\mathbb{D}|f'(z)|^2\,\omega(z)\,\operatorname{d}\!A(z)<\infty\}.\]
We say that $h(z)=\sum_n c_n z^n$ is a Hadamard multiplier for $\mathcal{D}_\omega$ if for every $f\in\mathcal{D}_\omega$ one has that $f\ast h\in\mathcal{D}_\omega$. 

In \cite{MR4032210} it was proven the following result:  given a sequence of complex numbers $(c_k)_k$, the function $h(z)=\sum_n c_n z^n$ is a Hadamard multiplier for $\mathcal{D}_\omega$ for every super-harmonic $\omega$, if and only if the infinite matrix
\begin{equation} \label{E:Definition of $T_c$}
 T_{c}=\left(\begin{array}{ccccc}
c_{0} & c_{1}-c_{0} & c_{2}-c_{1} & c_{3}-c_{2} & \ldots \\
0 & c_{1} & c_{2}-c_{1} & c_{3}-c_{2} & \cdots \\
0 & 0 & c_{2} & c_{3}-c_{2} & \ldots \\
0 & 0 & 0 & c_{3} & \ldots \\
\vdots & \vdots & \vdots & \vdots & \ddots
\end{array}\right)
\end{equation}
defines a bounded operator on $\ell^2$. In their work, Mashreghi and Ransford discussed sufficient conditions and necessary conditions for the boundedness of $T_c$, and certain estimates for its operator norm. In particular, their work led to a detailed study of approximation schemes in such spaces. They also left as an open problem to completely characterize the boundedness of the matrix $T_c$ on $\ell^2$, see Question 3 in \cite[Section 8]{MR4032210}.

Using the analysis of Rhaly matrices developed in Sections \ref{S:Section2} and \ref{S:Schatten}, we are able to fully answer the question of when the matrix $T_c$ is bounded. In fact, we completely characterize boundedness, compactness and Schatten class properties for $T_c$. The main idea is that, ignoring the diagonal entries, the resulting matrix can be studied as a Rhaly matrix. 

\begin{theorem}\label{T:main4}
Given $(c_k)_{k\in\mathbb{N}}$ a sequence of complex numbers, let $T_c$ be the infinite matrix defined in \eqref{E:Definition of $T_c$}. Then, the following properties hold.
\begin{enumerate}
    \item $T_c$ is bounded on $\ell^2$ if and only if
    \[\sup_{k\in\mathbb{N}} |c_k|<\infty \quad \text{and}  \quad \sup_{n\in\mathbb{N}} 2^n \sum_{j=2^n}^{2^{n+1}-1}|c_j-c_{j+1}|^2 <\infty.\]
    \item $T_c$ is compact on $\ell^2$ if and only if
    \[\lim_{k\to \infty} |c_k|=0 \quad \text{and}  \quad \lim_{n\to \infty} 2^n \sum_{j=2^n}^{2^{n+1}-1}|c_j-c_{j+1}|^2 =0.\]
    \item For $1<q<\infty$, $T_c$ belongs to $\mathcal{S}^q(\ell^2)$ if and only if
    \[\sum_{k\in\mathbb{N}} |c_k|^q<\infty \quad \text{and}  \quad \sum_{n\in\mathbb{N}} \biggl( 2^n \sum_{j=2^n}^{2^{n+1}-1}|c_j-c_{j+1}|^2 \biggr)^{\frac{q}{2}} <\infty.\]
\end{enumerate}
\end{theorem}

\section{Boundedness and compactness of $R_\alpha$} \label{S:Section2}
 From now on, we fix a sequence $\alpha=(\alpha_n)_n \in \ell^2$. All of the quantities that we introduce depend on such $\alpha$, but we will not stress it in the notation.

\begin{definition}
    For $n\in\mathbb{N}$, we define
\begin{equation*}    
J_{n}:=\sup_{m\geq n}\, (m+1-n)^{\frac{1}{2}}\biggl(\sum_{k=m}^{\infty}|\alpha_k|^2\biggr)^{\frac{1}{2}}.
\end{equation*}
Also, we denote $\sigma_{-1}:=|\alpha_0|$ and, for $k\in\mathbb{N}$,
\begin{equation}\label{D:defsigma}
    \sigma_{k}:=\biggl(\sum_{j=2^k}^{2^{k+1}-1}(j+1)|\alpha_j|^2\biggr)^{\frac{1}{2}}.
\end{equation}
\end{definition}

We note that, for $k\geq 0$,  
\begin{equation}\label{E:sigmabounds}
    2^{k} \sum_{j=2^k}^{2^{k+1}-1}|\alpha_j|^2\leq \sigma_k^2 \leq 2^{k+1} \sum_{j=2^k}^{2^{k+1}-1}|\alpha_j|^2.
\end{equation}

\subsection{Boundedness of $R_\alpha$ }
First of all, we recall the characterization of Bennett of the boundedness of \emph{factorable} matrices, that is, matrices that are defined as in Equation \eqref{factorable matrix} below.
\begin{theorem}[Theorem 2 of \cite{bennett}] \label{T:factorablebound}
    Let 
    \begin{equation}\label{factorable matrix}
    A=
\begin{pmatrix}
\alpha_0\beta_0 & 0 & 0 & 0 & 0 & \cdots\\
\alpha_1\beta_0 & \alpha_1\beta_1 & 0 & 0 & 0 & \cdots\\
\alpha_2\beta_0 & \alpha_2\beta_1 & \alpha_2\beta_2 & 0 & 0 & \cdots \\
\alpha_3\beta_0 & \alpha_3\beta_1 & \alpha_3\beta_2 & \alpha_3\beta_3 & 0 & \cdots\\
\vdots & \vdots & \vdots & \vdots &\vdots & \ddots
\end{pmatrix}
\end{equation}
    be a factorable matrix with $\alpha_k,\beta_k \in \mathbb{C}$, for $k\in\mathbb{N}$. The following conditions are equivalent.
    \begin{enumerate}[label=(\alph*)]
        \item The matrix $A\colon \ell^2\to\ell^2$ is bounded.
        \item It holds \begin{equation} \label{E:Gbddfactii}
            K_2 = \sup_{m\in\mathbb{N}} \left(\sum_{k=m}^\infty |\alpha_k|^2\right)^\frac{1}{2}\left(\sum_{j=0}^m |\beta_j|^{2}\right)^\frac{1}{2}<\infty.
        \end{equation}
    \end{enumerate}
    Moreover, we have also that
    \begin{equation}\label{E:Bennettbound}
      K_2 \leq \|A\| \leq 2\sqrt{2}K_2.
    \end{equation}
\end{theorem}

Theorem \ref{T:factorablebound} is used to establish the relationship between the boundedness of $R_\alpha$ and the quantities $J_n$, $\sigma_k$.
\begin{theorem}\label{T:boundgarnett}
The Rhaly operator $R_\alpha$ is bounded on $\ell^2$ if and only if $J_{0} < \infty$. In this case, 
 \[ J_0\leq \|R_{\alpha}\|\leq 2\sqrt{2} J_0.\]
\end{theorem}
\begin{proof}
We notice that $J_0$ is equal to $K_2$ of \eqref{E:Gbddfactii} in the case when $\beta_j=1$ for each $j$.
\end{proof}

\begin{theorem}\label{T:boundRasigma}
The Rhaly operator $R_\alpha$ is bounded on $\ell^2$ if and only if the sequence $(\sigma_k)_{k\in\mathbb{N}}$ is bounded. In particular,
\[
\frac{1}{\sqrt{2}}\sup_{h \geq -1}\sigma_h\leq \|R_\alpha\|\leq 4\sqrt{2} \sup_{h \geq -1}\sigma_h.
\]
\end{theorem}
  
\begin{proof}
Recall that \[J_0=\sup_{m\geq 0}\, (m+1)^\frac{1}{2}\biggl(\sum_{k=m}^\infty |\alpha_k|^2\biggr)^\frac{1}{2}.\]
If $m=0$, we have that
\[
\sum_{k=0}^\infty |\alpha_k|^2\leq |\alpha_0|^2+\sum_{h=0}^\infty\sum_{k=2^h}^{2^{h+1}-1}|\alpha_k|^2\leq |\alpha_0|^2+ \sum_{h=0}^\infty 2^{-h} \sigma_h^2\leq 3\sup_{h\geq -1}\sigma_h^2.
\]
Let $m\geq 1$, and let $k$ be the unique natural number such that $m\in[2^k,2^{k+1}-1].$ Then,
\begin{align*}
   (m+1)\sum_{j=m}^\infty |\alpha_j|^2 &\leq 2^{k+1}\sum_{j=2^k}^\infty |\alpha_j|^2 = 2^{k+1}\sum_{h=k}^\infty\sum_{j=2^h}^{2^{h+1}-1} |\alpha_j|^2\\
   &\leq  2^{k+1}\sum_{h=k}^\infty 2^{-h} \sigma_h^2\leq4\sup_{h\geq k} \sigma_h^2.
\end{align*}
It follows that $J_0\leq 2 \sup_{h \geq -1}\sigma_h$ and, by Theorem \ref{T:boundgarnett}, $R_\alpha$ is bounded by $\|R_\alpha\|\leq 4\sqrt{2}\sup_{h \geq -1}\sigma_h$. For the opposite inequality, we note that $\sigma_{-1}^2=|\alpha_0|^2\leq J_0^2$ and 
\[ \sigma_{k}^{2}\leq 2^{k+1}\sum_{j=2^k}^{2^{k+1}-1}|\alpha_j|^{2}\leq 2J_{0}^{2}.\]
Then, if $R_\alpha$ is bounded, $\sup_{h \geq -1}\sigma_h\leq \sqrt{2} J_0 \leq \sqrt{2} \|R_\alpha\|.$
\end{proof}

We are now ready to prove the first main theorem.

\begin{proof}[Proof of Theorem \ref{T:main1}]

The proof follows from Theorem \ref{T:boundRasigma}, the characterization of $\Lambda^{2}_{\frac{1}{2}}$ in terms of $\Delta_n$ and Equation \eqref{E:sigmabounds}.
\end{proof}

\subsection{Compactness of $R_\alpha$}
We start with a preliminary lemma. Although this result is well-known, we provide a proof for the convenience of the reader. From now on, we consider the elements of $\ell^2$ as square-summable functions $f\colon\mathbb{N}\to\mathbb{C}$. We say that a function $f\in\ell^2$ has \emph{finite support} if there exist natural numbers $a<b$ such that  $f(j)=0$ for every $j\notin I=[a,b]$. In this case, we say that $f$ is supported in $I$ or that $\operatorname{supp}(f)\subseteq I$. 

\begin{lemma}\label{L:approx}
    Let $P\in\mathcal{F}_N$ be a finite rank operator. Then, for every $\epsilon>0$ there exists an operator $Q$ with the following properties:
    \begin{itemize}
        \item $\|P-Q\|\leq\epsilon$;
        \item There exists $K\in\mathbb{N}$ such that, for every $f\in\ell^2$, \[\operatorname{supp}(Qf)\subseteq [0,K].\]
    \end{itemize}
\end{lemma}
\begin{proof}
We follow the argument of \cite[Lemma 2.2]{Evans1987}. Let $\{u_1,\ldots,u_N\}$ be elements of $\ell^2$ with the following property: for every $f\in\ell^2$ there exist unique coefficients $c_i(f),\,\,i=1,\ldots,N$ such that
\[Pf=\sum_{i=1}^N c_i(f)u_i.\]
Since the range of $P$ is finite-dimensional, every norm is equivalent, and there exists $M>0$ such that
\[\sum_{i=1}^N |c_i(f)|^2\leq M^2 \|Pf\|_2^2\leq M^2 \|P\|^2\|f\|_2^2.\]
By the density of the space of finitely supported sequences $c_0$ in $\ell^2$, there exist $\phi_1,\ldots,\phi_N\in c_0$ such that
\[\forall i=1,\ldots,N \qquad \|u_i-\phi_i\|_2\leq\frac{\epsilon}{M\|P\|\sqrt{N}}. \]
Set $Qf=\sum_{i=1}^{N}c_{i}(f)\phi_{i}$. Thus,
\begin{align*}
\|Pf-Qf\|_2^2 &\leq  \sum_{k\in\mathbb{N}} \biggl(\sum_{i=1}^N| c_i(f)(u_i(k)-\phi_i(k))|\biggr)^2\\
&\leq  N\sum_{i=1}^N| c_i(f)|^2 \|u_i-\phi_i\|_2^2\\
&\leq \epsilon^2 \|f\|_2^2.
\end{align*}
This proves that $\|P-Q\|\leq \epsilon$. The natural number $K$ can be chosen so that for every $j>K$ and every $1\leq i\leq N$, $\phi_i(j)=0.$
\end{proof}

We describe the compactness of $R_\alpha$ in terms of $J_n$. We note that $(J_n)_n$ is a decreasing sequence.

\begin{theorem}\label{Comp}
Let $\alpha=(\alpha_k)_k\in\ell^2$. Then, the inequalities
\[\lim_{n\to\infty}J_{n}\leq \|R_{\alpha}\|_e\leq 2\sqrt{2}\lim_{n\to\infty}J_{n}\]
hold. In particular, $R_\alpha$ is compact if and only if 
$\lim_{n\to\infty}J_n = 0.$
\end{theorem}
\begin{proof}
    Given any natural $N\geq 1$, we define the maps $R_{\alpha,N}$ and $P_{N}$ by
\begin{align}
\label{E:Truncated R}
    R_{\alpha,N}f(k)&=\alpha_k\chi_{[N,\infty)}(k)\sum_{j=0}^{k}f(j)\chi_{[N,\infty)}(j),\\
    \nonumber P_{N}f(k)&=\alpha_k\sum_{j=0}^{k}f(j)\chi_{[0,N)}(j),
\end{align}
where $\chi_{I}$ is the characteristic function of the interval $I$. We note that
\begin{align*}
R_\alpha f(k)&=\alpha_k\sum_{j=0}^{k}f(j)\\
&= \alpha_k\sum_{j=0}^{k}f(j)\chi_{[0,N)}(j)+\alpha_k\sum_{j=0}^{k}f(j)\chi_{[N,\infty)}(j)\\
&=\alpha_k\sum_{j=0}^{k}f(j)\chi_{[0,N)}(j) +\alpha_k\chi_{[N,\infty)}(k)\sum_{j=0}^{k}f(j)\chi_{[N,\infty)}(j)\\
&=P_Nf(k)+R_{\alpha,N}f(k).
\end{align*}
The operator $P_N$ is bounded. Indeed, for all $f\in \ell^{2}$,
\begin{align*}
 \|P_{N}f\|_{2}^2 &\leq \sum_{k=0}^{\infty}|\alpha_k|^{2}\biggl(\sum_{j=0}^{N-1}|f(j)|\biggr)^2\leq\|\alpha\|_{2}^2 N \|f\|_{2}^2.
\end{align*}
$P_{N}$ is also compact. We show that, if $(f_n)_n$ is a sequence in the closed unit ball of $\ell^2$, then $\{P_{N}f_{n}:n\in\mathbb{N}\}$ is relatively compact in $\ell^2$. Since $\{P_{N}f_{n}:n\in\mathbb{N}\}$ is bounded, it is enough to show that  
\[\forall\epsilon>0 \quad \exists N_\epsilon\in\mathbb{N} \,\,\,\text{such that}\,\,\, \forall n\in\mathbb{N} \quad \sum_{k=N_\epsilon}^\infty |(P_N f_n)(k)|^2 <\epsilon.\]
For every $\epsilon>0$, such $N_\epsilon$ can be easily obtained by the same argument used for the boundedness of $P_N$, and the fact that $\alpha\in\ell^2$.

As $R_{\alpha}=R_{\alpha,N} + P_{N}$, it follows that
\[ \|R_\alpha\|_e=\|R_{\alpha,N}\|_e.\]
From (\ref{E:Truncated R}) we see that
\[R_{\alpha,N}f(k)=
\begin{cases} 
      0, & 0\leq k< N,\\
      \alpha_k\sum_{j=N}^{k}f(j), &  k \geq N,
   \end{cases}
   \]
that is, $R_{\alpha,N}$ is the factorable matrix associated to $\alpha$ and $\beta^{(N)}$ defined as
\[\beta^{(N)}_j=\begin{cases} 0, &\text{if}\,\, j< N,\\
1, &\text{if}\,\, j\geq N.    
\end{cases}\]
By Theorem \ref{T:factorablebound},
\[\|R_\alpha\|_e=\|R_{\alpha,N}\|_e\leq \|R_{\alpha,N}\| \leq 2\sqrt{2} J_{N}.\]
Since $N$ may be chosen arbitrarily large, we see that 
\[\|R_{\alpha}\|_e\leq 2\sqrt{2}\lim_{n\to\infty}J_{n}.\]

We now establish the lower bound for $\|R_\alpha\|_e$. Let  $\lambda > \|R_{\alpha}\|_e$. Then, there exists $P\in \mathcal{F}$ such that for all $f\in\ell^2$,
\[\|R_{\alpha}f-Pf\|_{2}\leq \lambda \|f\|_{2}. \]
In light of Lemma \ref{L:approx}, we may assume that there exists $K\in\mathbb{N}$ such that, for all $f\in\ell^2$, $\operatorname{supp}Pf \subseteq [0,K-1]$. Hence,
\[\sum_{j=K}^{\infty}|R_{\alpha}f(j)|^{2}\leq\lambda^{2}\|f\|_2^2\]
for all $f\in\ell^{2}$. Now, let $a, b \in \mathbb{N}$ be such that $K< a < b$. Define $f\in\ell^{2}$ as the characteristic function $f=\chi_{[a,b]}$. Then,
\begin{align*}
\sum_{j=K}^{\infty}|R_{\alpha}f(j)|^{2}&\geq (b+1-a)^2  \sum_{j=b}^{\infty} |\alpha_j|^{2}.
\end{align*}
Since $\|f\|_{2}=(b+1-a)^{\frac{1}{2}}$, we have that
\[(b+1-a)^{\frac{1}{2}}\biggl(\sum_{j=b}^{\infty}|\alpha_k|^{2}\biggr)^{\frac{1}{2}}\leq\lambda.\]
Since $b>a$ may be chosen arbitrarily, we have that $\lim_{a\to\infty} J_a \leq \lambda$. Since $\lambda$ may be chosen arbitrarily close to $\|R_{\alpha}\|_e$, we conclude that $\lim_{n\to\infty} J_n \leq \|R_{\alpha}\|_e$. 
\end{proof}

\begin{theorem}\label{T:compRasigma}
   The map $R_\alpha$ is compact if and only if $\lim_{k\to \infty} \sigma_{k}= 0$. 
\end{theorem}
\begin{proof}
First, we show that for every $k\in\mathbb{N}$ fixed, $J_{2^k}\leq 4 \sup_{n\geq k}\sigma_n$. For $m\geq 2^k$, let $\ell\in\mathbb{N}$ be the unique natural number such that $2^\ell\leq m\leq2^{\ell+1}-1$. Notice that $\ell\geq k$. Then, as in the proof of Theorem \ref{T:boundRasigma},
\begin{align*}
   (m+1-2^k)\sum_{j=m}^\infty |\alpha_j|^2 &\leq (2^{\ell+1}-2^k)\sum_{j=2^\ell}^\infty |\alpha_j|^2\\
   &\leq 2^{\ell+1} \sum_{h=\ell}^\infty\sum_{j=2^h}^{2^{h+1}-1} |\alpha_j|^2\leq 4\sup_{n\geq k}\sigma_n^2.
\end{align*}
This shows that $J_{2^k}\leq 2\sup_{n\geq k}\sigma_n$. On the other hand, for $n\geq k+1$,

\begin{align*}
    \sigma_{n}^2&\leq \frac{2^{n+1}}{2^{n}+1-2^k} ( 2^{n}+1-2^k)\sum_{j=2^n}^{\infty}|\alpha_j|^2\\
    &\leq \frac{2^{n+1}}{2^n-2^{n-1}} ( 2^{n}+1-2^k)\sum_{j=2^n}^{\infty}|\alpha_j|^2\\
    &\leq 4 (2^{n}+1-2^k) \sum_{j=2^{n}}^\infty |\alpha_j|^2\leq 4 \sup_{m\geq 2^k} (m-2^k)\sum_{j=m}^\infty|\alpha_j|^2.
\end{align*}
This proves that for $n\geq k+1$ the inequality $\sigma_{n}\leq 2 J_{2^{k}}$ holds.
Hence $\lim_{a\to\infty} J_a = 0$ if and only if $\lim_{n\to\infty} \sigma_{n} = 0$, since $(J_a)_{a\in\mathbb{N}}$ is a decreasing sequence. The corollary now follows from Theorem \ref{Comp}.
\end{proof}

We are now ready to prove the second main theorem.

\begin{proof}[Proof of Theorem \ref{T:main2}]
With the help of Theorem \ref{T:compRasigma},  the result follows. 
\end{proof}

\section{Finer analysis on natural intervals}
In this section, we gather some technical tools to prove Theorem \ref{T:main3}. We adapt the analysis that Edmunds, Evans and Harris did on a class of integral operators on $\mathbb{R}$ in \cite{Edmunds1988, Edmunds1997}, to our discrete setting. Although this translation process might look literal, it is far from being straightforward.
In a nutshell, the continuum property of $\mathbb{R}$ plays a key role in their work. Passing to our discrete setting, many important properties get lost and non-trivial adjustments have to be made.

We consider intervals of natural numbers, that is, $I=[a,b]\cap\mathbb{N}$ with $a,b\in\mathbb{N}$. We denote the cardinality of such intervals as $\# I=b+1-a$.  By $\ell^{2}(I)$, we mean 
\[\ell^2(I)=\{f\in\ell^2\colon \operatorname{supp}(f)\subseteq I\}.\]
We denote the norm in $\ell^2(I)$ as 
\[\|f\|_{2,I}=\biggl(\sum_{k \in I}|f(k)|^2\biggr)^{\frac{1}{2}}.\]
It is convenient to consider the measure $\mu=\mu_\alpha$ defined as
\[\mu(I):=\sum_{k\in I}|\alpha_k|^2, \qquad  I\subseteq \mathbb{N}.\]
Finally, we also introduce the notation
\[\|f\|_{2,I,\mu}=\biggl(\sum_{k \in I}|f(k)|^2|\alpha_k|^2\biggr)^{\frac{1}{2}}.\]
Given  $I\subseteq \mathbb{N}$ a natural interval, for any $f \in \ell^2(\mathbb{N})$, we denote 
\begin{align*}
    l(I,f)&:=\sum_{k\in I}\sum_{n\in I\setminus\{k\}}\biggl|\alpha_k \alpha_n\sum^{\max(k,n)}_{j=\min(k,n)+1}f(j)\biggr|^{2}.
\end{align*}

\begin{definition}
Let $I=[a,b]$ be an interval in $\mathbb{N}$. We define
\begin{equation}\label{L(I) def} 
    L(I):=\biggl(\sup_{\|f\|_{2,I}\leq 1}\left\lbrace\frac{l(I,f)}{\sum_{k\in I}|\alpha_k|^2}\right\rbrace \biggr)^{\frac{1}{2}}=\biggl(\sup_{\|f\|_{2,I}\leq 1}\left\lbrace\frac{l(I,f)}{\mu(I)}\right\rbrace \biggr)^{\frac{1}{2}}.
\end{equation}
If $\alpha_k=0$ for every $k\in I$, we set $L(I)=0$. 
\end{definition}
Note that, in the case where $a=b$ and the interval $I$ is just the singleton $\{a\}$, then $l(I,f)=0$ for every $f$ and $L(I)=0$.

\begin{lemma}\label{lemma 5}
If $I=[a,b]$ with $a,b\in\mathbb{N}$, $a<b$, then the quantity $L([a, b])$ decreases as $a$ increases and it increases as $b$ increases.
\end{lemma}
\proof
Notice that for each $f$ we have that $l(I,f)\leq l(I,|f|)$. When evaluating the supremum $L(I)$, we may consider without loss of generality only positive functions. Using Fubini's theorem, we have that
\begin{align*}
     l(I,f)&=\sum_{k\in I}\sum_{n\in I\setminus\{k\}}\biggl|\alpha_k \alpha_n\sum^{\max(k,n)}_{j=\min(k,n)+1}f(j)\biggr|^{2}\\
    &=\sum_{k=a}^{b}\sum_{n=a}^{k-1}\biggl|\alpha_k \alpha_n\sum^k_{j=n+1}f(j)\biggr|^{2}+\sum_{k=a}^{b}\sum_{n=k+1}^b\biggl|\alpha_k\alpha_n\sum^n_{j=k+1}f(j)\biggr|^{2}\\
    &=\sum_{n=a}^{b-1}\sum_{k=n+1}^{b}\biggl|\alpha_k \alpha_n\sum^k_{j=n+1}f(j)\biggr|^{2}+\sum_{k=a}^{b}\sum_{n=k+1}^b\biggl|\alpha_k\alpha_n\sum^n_{j=k+1}f(j)\biggr|^{2}\\
&=2\sum_{k=a}^{b}\sum_{n=k+1}^b\biggl|\alpha_k\alpha_n\sum^n_{j=k+1}f(j)\biggr|^{2}.
\end{align*}
We introduce the function
\[g(k)=\sum_{n=k+1}^b\biggl|\alpha_n\sum^n_{j=k+1}f(j)\biggr|^2.\]
Since $f\geq 0$, the function $g$ is decreasing. Now, let $a_1\in ( a,b]$.
\begin{align*}
 \frac{l([a,b],f)}{\mu([a,b])}&=\frac{2}{\mu([a,b])} \sum_{k=a}^{b}|\alpha_k|^2 g(k)\\
 &\geq \frac{2}{\mu([a,b])}   \biggl(g(a_1-1)\sum_{k=a}^{a_1-1}|\alpha_k|^2+\sum_{k=a_1}^{b}|\alpha_k|^2g(k)
 \biggr).
\end{align*}

Notice that, again by the monotonicity of $g$,
\[\frac{1}{\sum_{k=a_1}^b |\alpha_k|^2}\sum_{k=a_1}^b g(k) |\alpha_k|^2 \leq g(a_1)\leq g(a_1-1).\]
It follows that
\begin{align*}
 \frac{l([a,b],f)}{\sum_{k=a}^b|\alpha_k|^2}&\geq \frac{2}{\sum_{k=a}^b|\alpha_k|^2} \biggl(\sum_{k=a_1}^{b}|\alpha_k|^2g(k) \biggr)\biggl(\frac{\sum_{k=a}^{a_1-1} |\alpha_k|^2}{\sum_{k=a_1}^b |\alpha_k|^2}+1\biggr) \\
 &= \frac{2}{\sum_{k=a_1}^b|\alpha_k|^2} \biggl(\sum_{k=a_1}^{b}|\alpha_k|^2g(k) \biggr)=\frac{l([a_1,b],f)}{\sum_{k=a_1}^b|\alpha_k|^2}.
\end{align*}
Hence $L([a, b])$ decreases as $a$ increases. In a similar way, we show that $L([a, b])$ increases as $b$ increases, writing
\begin{equation*} \label{E:errore1}
    l([a,b],f)=2\sum_{k=a}^{b}\sum_{n=a}^{k-1}\biggl|\alpha_k\alpha_n\sum^k_{j=n+1}f(j)\biggr|^{2}.\qedhere
\end{equation*}

To deal with infinite intervals we set
\[L([a, \infty)) := \lim_{b\to\infty} L([a, b])=\sup_{b\geq a} L([a,b]). \]

\begin{definition}\label{D:JAcBc}
Let $ I=[a,b]$ be an interval in $\mathbb{N}$. Then 
\begin{equation}\label{J(I) def}
  J(I):=\inf_{c \in I} \lbrace  \max (A_{c}^I,B_{c}^I)\rbrace ,
\end{equation}
where
\begin{equation*}\label{A(c) def}
  A_{c}^I:=\sup_{a\leq s\leq c}(c-s)^{\frac{1}{2}}\biggl(\sum_{k=a}^{s}|\alpha_k|^{2}\biggr)^{\frac{1}{2}}, \quad\,  B_{c}^I:=\sup_{c\leq s\leq b}(s-c)^{\frac{1}{2}}\biggl(\sum_{k=s}^{b}|\alpha_k|^{2}\biggr)^{\frac{1}{2}}.
\end{equation*}
\end{definition}
We will omit $I$ from the notation $A_c^I,B_c^I$ when it is clear. We note that $A_{a}^I=0=B_{b}^I$.

\begin{notation}
Given a function $f\in \ell^2$ and a natural interval $I$ with $\mu(I)>0$, we denote
\[f_{I}:=\frac{1}{\mu(I)}\int_I f\operatorname{d}\!\mu=\frac{1}{\sum_{k\in I}|\alpha_k|^2}\sum_{k\in I}f(k)|\alpha_k|^{2}\]
the average over $I$ of $f$ with respect to the measure $\mu$. Also, for every $k\in\mathbb{N}$ we write the cumulative function
\[S_f(k)=\sum_{j=0}^{k}f(j).\]
\end{notation}
Notice that $R_\alpha f(k)=\alpha_k S_f(k).$

\begin{definition}
    For a non-empty interval $I$, we write
\begin{equation}
K(I) := \sup_{\|f\|_{2,I}\leq 1} \biggl(\sum_{k\in I} |S_f(k)-(S_f)_I|^2 |\alpha_k|^2\biggr)^{\frac{1}{2}}.
\end{equation}
\end{definition}
The next lemma is crucial, as it establishes relations between the quantities $L(I), K(I), J(I)$. 
\begin{lemma}\label{L:relationKJL}
Let $I=[a,b]$ be an interval in $\mathbb{N}$. Then,   
\begin{equation}\label{K and L equi}
    \frac{1}{2}J(I)\leq K(I) = \frac{1}{\sqrt{2}}L(I)\leq 2J(I).
\end{equation}
\end{lemma}
\begin{proof}
Let $f$ with $\|f\|_{2,I}\leq 1$. Setting $G=S_f-(S_f)_I$, one has that the average  $G_I=0$, and
\begin{equation*}\label{eq L and K}
\begin{split}
   \int_{I}\int_{I}|S_f(k)-S_f(n)|^{2}\,\operatorname{d}\!\mu(k)\operatorname{d}\!\mu(n)&=\int_{I}\int_{I}|G(k)-G(n)|^{2}\,\operatorname{d}\!\mu(k)\operatorname{d}\!\mu(n)\\
   &=2\mu(I)\int_I |G(k)|^2\operatorname{d}\!\mu(k)\\
    &=2\mu(I)\|S_f-(S_f)_{I}\|^{2}_{2,I,\mu}.
    \end{split}
\end{equation*}
Thus, 
\begin{align*}
\|S_f-(S_f)_{I}\|^{2}_{2,I,\mu}&=\frac{1}{2\mu(I)} \sum_{n\in I}\sum_{k \in I\setminus\{n\}}|S_f(k)-S_f(n)|^{2}|\alpha_k|^{2}|\alpha_n|^{2} \\
&=\frac{1}{2\mu(I)}\sum_{n\in I}\sum_{k\in I\setminus\{n\}}\biggl|\sum_{j=\min(n,k)+1}^{\max(n,k)}f(j)\biggr|^{2}|\alpha_k|^{2}|\alpha_n|^{2}\\
&=\frac{1}{2} \frac{l(I,f)}{\mu(I)}.
\end{align*}
Taking now the supremum over $\|f\|_{2,I}\leq 1$ in the left-hand side, we conclude that $K(I)=L(I)/ \sqrt{2}$.

We move on to the lower bound for $K(I)$. We have that
\begin{align*}
\mu(I)\|S_f-(S_f)_{I}\|^{2}_{2,I,\mu}&=\frac{1}{2} \sum_{k \in I}\sum_{l \in I}|S_f(k)-S_f(l)|^{2}|\alpha_k|^{2}|\alpha_l|^{2} \\
&=\sum_{l \in I}\sum_{k< l}|S_f(k)-S_f(l)|^{2}|\alpha_k|^{2}|\alpha_l|^{2} \\
&=\sum_{l \in I}\sum_{k<l}\biggl|\sum_{j=k+1}^{l}f(j)\biggr|^{2}|\alpha_k|^{2}|\alpha_l|^{2}.
\end{align*}

We consider the set 
\[E=\{n\in (a,b] \colon \sum_{j=a}^{n-1}|\alpha_j|^2 < \frac{1}{2}\sum_{j=a}^b|\alpha_j|^2\}.\]
If $E$ is non-empty, we set $N=\max E$. Otherwise, we set $N=a$.
The number $N$ has the following properties:
\begin{align}
    \label{E:relationKJLN1}\sum_{j=a}^{N}|\alpha_j|^2 &\geq \frac{1}{2}\sum_{j=a}^b|\alpha_j|^2 ,\\
   \label{E:relationKJLN2} \sum_{j=N}^{b} |\alpha_j|^2&=\sum_{j=a}^{b} |\alpha_j|^2-\sum_{j=a}^{N-1} |\alpha_j|^2 \geq \frac{1}{2}\sum_{j=a}^{b} |\alpha_j|^2.
\end{align}
Assume that $N>a$. For $a\leq K< N$, consider the characteristic function $g=\chi_{[K,N]}$. Then,
\begin{align*}
    \mu(I)\|S_g-(S_g)_{I}\|^{2}_{2,I,\mu}&=\sum_{l \in I}\sum_{i<l}\biggl|\sum_{j=i+1}^{l}\chi_{[K,N]} (j)\biggr|^{2}|\alpha_i|^{2}|\alpha_l|^{2}\\
     & \geq \sum_{l \geq N} \sum_{i\leq K}\biggl|\sum_{j=i+1}^{l}\chi_{[K,N]} (j)\biggr|^{2}|\alpha_i|^{2}|\alpha_l|^{2}\\
    & \geq \sum_{l \geq N} \sum_{i\leq K}(N-K)^{2}|\alpha_i|^{2}|\alpha_l|^{2}.
\end{align*}
Therefore, by \eqref{E:relationKJLN2},
\begin{align*}
\|S_g-(S_g)_{I}\|_{2,I,\mu}^2&\geq (N-K)^2 \biggl(\sum_{i=a}^{K}|\alpha_i|^{2}\biggr)\biggl( \frac{\sum_{l=N}^b |\alpha_l|^{2}}{\sum_{l=a}^{b}|\alpha_l|^{2}}\biggr)\\
    &\geq\frac{1}{2} (N-K)^2 \biggl(\sum_{i=a}^{K}|\alpha_i|^{2}\biggr),
\end{align*}
and it follows that
\begin{align*}
    K(I)^2 &\geq \frac{1}{\|g\|_{2,I}^2} \|S_g-(S_g)_{I}\|_{2,I,\mu}^2\geq \frac{N-K}{2(N-K+1)} (N-K) \biggl(\sum_{i=a}^{K}|\alpha_i|^{2}\biggr)\\
    &\geq \frac{1}{4} (N-K) \biggl(\sum_{i=a}^{K}|\alpha_i|^{2}\biggr),
\end{align*}
 uniformly on $K \in [a,N)$. We conclude that
\[A_N^I =\sup_{a\leq K\leq N}(N-K)^\frac{1}{2} \biggl(\sum_{i=a}^{K}|\alpha_i|^{2}\biggr)^{\frac{1}{2}}\leq 2 K(I).\]
If $N=a$, then $A_N=0$. Analogously, for $b\geq M> N$, choosing $h=\chi_{[N,M]}$ we have that
\begin{align*}
    \mu(I)\|S_h-(S_h)_{I}\|^{2}_{2,I,\mu}&\geq\sum_{l \geq M} \sum_{i\leq N}(M-N)^{2}|\alpha_i|^{2}|\alpha_l|^{2},
\end{align*}
and we obtain that
\begin{align*}
    \frac{\|S_h-(S_h)_{I}\|^{2}_{2,I,\mu}}{\|h\|_{2,I}^2} &\geq \frac{(M-N)^2}{(M-N+1)} \biggl(\sum_{l=M}^b|\alpha_l|^2\biggr) \frac{\sum_{i=a}^N |\alpha_i|^2}{\sum_{i=a}^b |\alpha_i|^2} \\
    &\geq \frac{M-N}{4}\biggl(\sum_{l=M}^b|\alpha_l|^2\biggr),
\end{align*}
where in the last inequality we used \eqref{E:relationKJLN1}. It follows that
\[B_N^I =\sup_{N\leq M\leq b}(M-N)^{\frac{1}{2}}\biggl(\sum_{l=M}^{b}|\alpha_l|^{2}\biggr)^{\frac{1}{2}}\leq 2 K(I).\]
In conclusion,
\[J(I)=\inf_{c\in[a,b]}\big(\max\{A_c,B_c\}\big)\leq \max\{A_N,B_N\} \leq 2 K(I).\]
We prove the upper bound for $K(I)$. This is based on an argument of Maz'ja \cite{Mazja}. Let $c\in \mathbb{N}$ be any point in $[a, b]$. Then,
\begin{equation}
    \label{first est}
       \sum_{k\in[a,c]} |S_f(k)-S_f(c)|^2 |\alpha_k|^2 
        \leq\sum_{k=a}^{c-1}|\alpha_k|^2\biggl(\sum_{l=k+1}^{c}|f_l| \biggr)^{2}.
\end{equation} 
Consider
\[g(j):=(c+1-j)^{1/4},\]
and observe that by the Cauchy-Schwarz inequality,
\begin{align*}
\biggl(\sum_{j=k+1}^{c}|f(k)|\biggr)^{2}
&\leq\biggl(\sum_{j=k+1}^{c}|f(j)|^{2}|g(j)|^{2}\biggr)\biggl(\sum_{j=k+1}^{c}( c+1-j)^{-\frac{1}{2}}\biggr)\\
&\leq\biggl(\sum_{j=k+1}^{c}|f(j)|^{2}|g(j)|^{2}\biggr)\biggl(\int_0^{c-k}x^{-\frac{1}{2}}\operatorname{d}\!x\biggr)\\
&=2\biggl(\sum_{j=k+1}^{c}|f(j)|^{2}|g(j)|^{2}\biggr)(c-k)^{\frac{1}{2}}.
\end{align*}

The use of this in (\ref{first est}) shows that
\begin{equation*}
\begin{split}
    \sum_{k\in[a,c]} |S_f(k)-S_f(c)|^2 |\alpha_k|^2 &\leq 2 \sum_{k=a}^{c-1}|\alpha_k|^{2}\biggl(\sum_{j=k+1}^{c}|f(j)|^{2}|g(j)|^{2}\biggr)(c-k)^{\frac{1}{2}} \\
    &= 2 \sum_{j=a+1}^{c}|f(j)|^{2}|g(j)|^{2}\sum_{k=a}^{j-1}|\alpha_k|^{2}(c-k)^{\frac{1}{2}},
\end{split}
\end{equation*}
where in the last equality we used Fubini's theorem. By the definition of $A_c$, that is,
\[ A_{c}=\sup_{a\leq k\leq c}(c-k)^{\frac{1}{2}}\biggl(\sum_{l=a}^{k}|\alpha_l|^{2}\biggr)^{\frac{1}{2}},\]
it follows that
\begin{align*}
 \sum_{k\in[a,c]} |S_f(k)-S_f(c)|^2 |\alpha_k|^2 &\leq 2A_c \sum_{j=a+1}^{c}|f(j)|^{2}|g(j)|^{2}\sum_{k=a}^{j-1}|\alpha_k|^{2}\biggl(\sum_{l=a}^{k}|\alpha_l|^{2}\biggr)^{-\frac{1}{2}}\\
&\leq 4 A_{c}\sum_{j=a+1}^{c}|f(j)|^2|g(j)|^{2}\biggl(\sum_{k=a}^{j-1}|\alpha_k|^{2}\biggr)^{\frac{1}{2}}\\
&= 4 A_{c}\sum_{j=a+1}^{c}|f(j)|^2(c-(j-1))^{\frac{1}{2}}\biggl(\sum_{k=a}^{j-1}|\alpha_k|^{2}\biggr)^{\frac{1}{2}} \\
&\leq 4A_{c}^{2}\sum_{j=a+1}^{c}|f(j)|^2,
\end{align*}
where in the second inequality we have used Bennett's power rule \cite[Lemma 1]{Edmunds1988} with $p=\frac{1}{2}$. Similarly, we find that
\begin{equation*}\label{estim with bc}
\begin{split}
   \sum_{k\in[c,b]} |S_f(k)-S_f(c)|^2 |\alpha_k|^2 &= \sum_{k\in(c,b]} |S_f(k)-S_f(c)|^2 |\alpha_k|^2\\
   &\leq 4B_{c}^2\sum_{k=c+1}^{b}|f(k)|^{2}.
\end{split}
\end{equation*}
Hence, for every $c\in I,$
\begin{align*}
    \sum_{k\in I} |S_f(k)-S_f(c)|^2 |\alpha_k|^2=& \sum_{k=a}^c |S_f(k)-S_f(c)|^2 |\alpha_k|^2 +\\
    &+\sum_{k=c+1}^b |S_f(k)-S_f(c)|^2 |\alpha_k|^2  \\
    \leq& 4\max(A_c,B_c)^{2}\|f\|^{2}_{2,I}.
\end{align*}
Using the minimizing property of the average
\[\|S_f-(S_f)_{I}\|_{2,[a,b],\mu}=\inf_{k \in \mathbb{C}}\|S_f-k\|_{2,[a,b],\mu},\]
we conclude that for every $f\in\ell^{2}(I)\setminus\{0\}$
\[\frac{1}{\|f\|_{2,I}} \|S_f-(S_f)_{I}\|_{2,I,\mu}\leq 2\inf_{c\in I}\max(A_c,B_c)=2 J(I).\]
This gives the second inequality in the statement of the lemma.
\end{proof}

\section{$(\epsilon,L)$-sequences} \label{S:eLseq}
We introduce the $(\epsilon,L)$-sequences, a key tool developed in \cite{Edmunds1988,Edmunds1997}.  
\begin{definition} 
Given a real number $\epsilon>0$, we choose natural numbers $c_k$ with the following rule: we set $c_0=0$, and for $k>0$
\begin{equation}\label{def of c_k}
     c_{k+1}=\inf\{t\in\mathbb{N}\colon t>c_{k}, \quad L([c_{k},t-1])>\epsilon\}.
\end{equation}
We say that these numbers form an
$(\epsilon, L)$-sequence. 
\end{definition}

If there exists a natural number $N$ such that $c_N\in\mathbb{N}$ and $c_{N+1}=+\infty$, we say that the $(\epsilon, L)$-sequence $\{c_k\}_{k=0}^N$ is finite and it has length $N$. Otherwise, we say that the $(\epsilon, L)$-sequence $\{c_k\}_{k\in\mathbb{N}}$ is infinite.

We recall that for every $a\in\mathbb{N}$, $L(\{a\})=0$. In particular, $c_{1}\geq 2$. The numbers $c_k$ forming an $(\epsilon,L)$-sequence have the following crucial properties:  for every $k\in\mathbb{N}$, we have that $c_{k+1}>c_{k}+1$, that $L([c_k,c_{k+1}-2])\leq \varepsilon$ and, if $c_{k+1}\in\mathbb{N}$, then $L([c_k,c_{k+1}-1])>\varepsilon$. 

We point out that in the continuous case studied in \cite{Edmunds1988,Edmunds1997}, the numbers forming an $(\epsilon,L)$-sequence are chosen to satisfy the stronger property $L([c_k,c_{k+1}])=\epsilon$ and $\{[c_k,c_{k+1}]\}_{k\in\mathbb{N}}$ is a partition of the real line, up to a set of Lebesgue measure $0$.  Due to the discrete nature of $\mathbb{N}$, in our context these properties are not satisfied.  

We relate the approximation numbers of $R_\alpha$ with the length of the $(\epsilon,L)$-sequences. The key idea is the following:  an $(\epsilon,L)$-sequence of length $N$ provides control from above over the approximation number $a_{2N+2}$, in terms of $\epsilon$, and from below over the approximation number $a_{N}$.

\begin{lemma}\label{L:singvaluppereL}
Let $\epsilon > 0$, and suppose that the associated $(\epsilon,L)$-sequence has finite length $N$. Then, we have that
\[ a_{2N+2}(R_\alpha)\leq \frac{\epsilon}{\sqrt{2}}.\]
\end{lemma}
\begin{proof}
By definition, there exist $N\in\mathbb{N}$ and numbers $c_k\in\mathbb{N}$, $k = 0, 1,\ldots, N$, with $0 = c_0 < c_{1} < \ldots < c_N < c_{N+1} =\infty$, such that $L(J_k)\leq \epsilon$
for $k = 0,1, ..., N$, where $J_k = [c_k, c_{k+1}-2]$\footnote{We point out that $J_N=[c_N, +\infty)$.}. Let $f\in \ell^{2}(\mathbb{N})$ be such that $\|f\|_2 = 1$, and write 
\[Pf:=\sum_{k=0}^{N}P_{J_{k}}f\]
where 
\[P_{J_k}f(j):= (S_f)_{J_k}\alpha_j \chi_{J_k}(j) +   (R_\alpha f)(c_{k+1}-1) \chi_{\{c_{k+1}-1\}}(j),\]
where as usual $\chi$ denotes a characteristic function, and \[S_f(j)=\sum_{m=0}^j f(m),\quad (S_f)_{J_k}=\frac{1}{\sum_{j\in J_k}|\alpha_j|^2}\sum_{j\in J_k}S_f(j)|\alpha_j|^{2}. \]
We interpret $\chi_{\{c_{N+1}-1\}}=\chi_{\{\infty\}}$ to be equal to $0$. Then, $P$ is a bounded linear operator on $\ell^{2}(\mathbb{N})$ and its range has dimension less than or equal to $2N+1$. Since $\{[c_k,c_{k+1}-1]\}_{k=1,\ldots,N}$ is a partition of $\mathbb{N}$ and each $P_{J_k}$ is supported in $[c_k,c_{k+1}-1]$, we have that
\begin{align*}
\|R_{\alpha}f-Pf\|^{2}_{2}&=\sum_{k=0}^{N}\sum_{j=c_k}^{c_{k+1}-1}|R_\alpha f(j) - P_{J_k}f(j)|^2\\
&=\sum_{k=0}^{N}\sum_{j=c_k}^{c_{k+1}-2}|R_\alpha f(j) - \alpha_j (S_f)_{J_k}|^2\\
&=\sum_{k=0}^{N}  \|S_f-(S_f)_{J_k}\|_{2,J_k,\mu}^2\\
&\leq\sum_{k=0}^{N} K^{2}(J_k)\|f\|^{2}_{2,J_{k}}.
\end{align*}
Note that we have used the identity $R_\alpha f(j)=\alpha_j S_f(j)$. Since $K(J_k)= \dfrac{L(J_{k})}{\sqrt{2}}\leq \dfrac{\epsilon}{\sqrt{2}}$ for Lemma \ref{L:relationKJL}, we have that
\[\|R_{\alpha}f-Pf\|^{2}_{2}\leq \frac{\epsilon^{2}}{2}\sum_{k=0}^{N}\|f\|^{2}_{2,J_{k}}\leq \frac{\epsilon^{2}}{2}\|f\|_2^2.\]
It follows that 
$a_{2N+2}(R_{\alpha})\leq \dfrac{\epsilon}{\sqrt{2}}. $
\end{proof}

\begin{lemma}\label{L:lowerboundaN}
Let $\epsilon > 0$ and suppose that the $(\epsilon,L)$-sequence has length at least $N\in\mathbb{N}$. Then, we have that
\[ a_{N}(R_{\alpha})\geq\frac{\epsilon}{\sqrt{2}}.\]
\end{lemma}
\begin{proof}
By definition, there exist $N\in\mathbb{N}$ and numbers $c_k\in\mathbb{N}$, $k = 1,\ldots, N$, with $0 = c_0 < c_{1} < \ldots < c_N <\infty$, such that $L(I_k)>\epsilon$
for $k = 0,1, ..., N-1$, where $I_k = [c_k, c_{k+1}-1]$. 

Let $\eta\in(0, 1)$. Since $K(I_k)=L(I_k)/\sqrt{2}$ and $L(I_{k})\geq \epsilon$, we see that for
$k = 0, 1, ..., N - 1 $ there exists $\phi_{k}\in\ell^{2}(I_{k})$ such that
\[\frac{\|S_{\phi_{k}}-(S_{\phi_{k}})_{I_{k}} \|_{2,I_{k},\mu}}{\|\phi_{k}\|_{2,I_{k}}}>\frac{\epsilon}{\sqrt{2}}\eta.\]
Set $\phi_{k}(j)= 0$ for all $j\in\mathbb{N}\setminus I_{k}$. Let $P\in \mathcal{F}_{N-1}$. Then, there exist $\lambda_0,\ldots,\lambda_{N-1}\in\mathbb{C}$ not all zero, such that
\[\sum_{k=0}^{N-1}\lambda_{k}P\phi_{k}=0.\]
Introducing
\[\phi:=\sum_{k=0}^{N-1}\lambda_{k}\phi_{k},\]
then 
\[S_\phi(m)=\sum_{j=0}^{m}\phi(j)=\sum_{j=0}^{m} \sum_{k=0}^{N-1}\lambda_{k}\phi_{k}(j),\qquad m\in\mathbb{N}.\]
Now, we fix $h\in \{0,\ldots,N-1\}$ and an interval $I_h=[c_h,c_{h+1}-1]$. If $m\in I_h$, then, for $j\leq m$ and every $k>h$, $\phi_k(j)=0$. Therefore,
\begin{align*}
S_\phi(m)=\lambda_h S_{\phi_{h}}(m) +\sum_{j=0}^{m}\sum_{k=0}^{h-1}\lambda_{k}\phi_{k}(j)=\lambda_hS_{\phi_{h}}(m)+\sum_{j=0}^{c_h}\sum_{k=0}^{h-1}\lambda_{k}\phi_{k}(j).
\end{align*}
In particular, for every $m\in I_h$ we can write $S_\phi(m)=\lambda_{h}S_{\phi_h}(m)+\mu_{h}$, for some $\mu_{h}$ that depends only on $h$ and not on $m$. 
Recall that 
\[R_\alpha \phi(m)=\alpha_m \sum_{j=0}^m\phi(j)=\alpha_m S_\phi(m).\]
Then, since the intervals $\{I_k\}_{k=0,\ldots N-1}$ are pairwise disjoint,
\[\| R_{\alpha}\phi-P\phi\|^{2}_{2}=\|R_{\alpha}\phi\|_{2}^{2}\geq \sum_{h=0}^{N-1}\|S_\phi\|^{2}_{2,I_{h},\mu}=\sum_{h=0}^{N-1}\|\lambda_{h}S_{\phi_h}+\mu_{h}\|^{2}_{2,I_{h},\mu}.\]
Using again the minimizing property of the average, we obtain that
\begin{align*}
 \| R_{\alpha}\phi-P\phi\|^{2}_{2}&\geq\sum_{h=0}^{N-1}\|\lambda_{h}S_{\phi_{h}}-(\lambda_{h}S_{\phi_{h}})_{I_{h}}\|^{2}_{2,I_{h},\mu}\\
 &\geq \frac{\eta^2\epsilon^2}{2}\sum_{h=0}^{N-1}|\lambda_{h}|^{2}\|\phi_{h}\|^{2}_{2,I_{h}}=\frac{\eta^2\epsilon^2}{2}\|\phi\|^{2}_{2}.
\end{align*}
Hence,
\[\|R_{\alpha}\phi-P\phi\|_{2}\geq \frac{\eta\epsilon}{\sqrt{2}}\|\phi\|_{2}.\]
which shows that $a_{N}(R_{\alpha})\geq \eta\epsilon /\sqrt{2}$. Since $\eta$ may be chosen arbitrarily close
to $1$, we conclude that $a_{N}(R_{\alpha})\geq \epsilon/\sqrt{2}$.
\end{proof}

The existence of $(\epsilon,L)$-sequences of finite length may be related to the boundedness and the compactness of the operator $R_\alpha$.

\begin{theorem}\label{T:boundRaeLseq} The following conditions are equivalent. 
\begin{enumerate}[label=(\roman*)]
    \item The operator $R_\alpha$ is bounded on $\ell^2$.
    \item $L([0, \infty)) < \infty$.
    \item There exists an $(\epsilon,L)$-sequence of length $N=0$.
\end{enumerate} 
\end{theorem}
\begin{proof}
We show that $(ii)\implies(iii)\implies(i)$. Take $\epsilon \geq L(0, \infty)$. In view of the monotonicity properties of $L$, we see that the associated $(\epsilon,L)$-sequence is given by $c_0=0,  c_1=\infty$. In particular, the length is $0$. Having now an $(\epsilon,L)$-sequence of length $0$, by Lemma \ref{L:singvaluppereL} the inequality
\[a_{2}(R_{\alpha})\leq \frac{\epsilon}{\sqrt{2}}\]
holds, that is,
\[ 
\inf_{P \in \mathcal{F}_1}\|R_{\alpha}-P\|\leq \frac{\epsilon}{\sqrt{2}}. 
\] 
Since every such operator $P$ is bounded, it follows that $R_\alpha$ is bounded as well.

To conclude, we prove that $(i)$ implies $(ii)$. We claim that, if $R_\alpha$ is bounded, then $L(I)\leq \sqrt{2}\|R_{\alpha}\|$ for any interval $I$ of the form $I=[0,b]$, $b\in\mathbb{N}$. Arguing by contradiction, suppose that $L([0,b])> \sqrt{2}\|R_{\alpha}\|+\delta$ for some natural number $b$ and $\delta>0$. Let 
$$
c_1:=\min\{b\in\mathbb{N}\setminus\{0\}\colon L([0,b-1])>\sqrt{2}\|R_{\alpha}\|+\delta\}<\infty.
$$ 
Now notice that $0=c_0< c_1$ is the start of a (possibly infinite) $(\epsilon,L)$-sequence with length at least $1$, where $\epsilon=\sqrt{2}\|R_\alpha\|+\delta$. Thus, by Lemma \ref{L:lowerboundaN}, 
\[
\|R_\alpha\|=a_1(R_\alpha)\geq\frac{\epsilon}{\sqrt{2}}>\|R_\alpha\|,
\]
which is impossible.
Then, we showed that $L([0,b])$ is uniformly bounded in $b$, and we conclude that $L([0, \infty)) \leq \sqrt{2} \|R_\alpha\| < \infty$.  
\end{proof}

By the monotonicity of $L$ (see Lemma \ref{lemma 5}), the previous theorem can be stated also in the following way: $R_\alpha$ is bounded if and only if $\sup_{n\in\mathbb{N}} L([n,\infty))<\infty.$ Also, similarly one can show that $R_\alpha$ is bounded if and only if $L([N,\infty))<\infty$ for some $N\in\mathbb{N}$.

\begin{theorem}\label{T:compRaeLseq1}
Suppose that $R_\alpha$ is bounded. Let $\mathcal{L}:= \lim_{n\to \infty} L([n, \infty)).$ Then, $\|R_{\alpha}\|_e= \dfrac{\mathcal{L}}{\sqrt{2}}$. In particular, $R_\alpha$ is compact if and only if $\mathcal{L}=0$.
\end{theorem}

\begin{proof}
First, notice that by monotonicity and Theorem \ref{T:boundRaeLseq}, we have that $\mathcal{L}\leq L(0,+\infty)<\infty$. We show that for every $\epsilon >\mathcal{L}$, the associated $(\epsilon,L)$-sequence has finite length. Because of the definition, there exists $A\in\mathbb{N}$ such that $L([A,+\infty]) \leq \epsilon$, and without loss of generality we also assume that $A$ is the smallest natural number with this property. Now, we define $M$ as
\[M:=\sup\{k\in\mathbb{N}\colon c_{k}\leq A\}<\infty.\]
The length of the $(\epsilon,L)$-sequence is now $N\leq M+1$. This is because $c_{M+1}$ may still be finite, but then necessarily
\[c_{M+2}=\inf\{t>c_{M+1}\colon L([c_{M+1},t-1])>\epsilon\}=+\infty,\]
since for every $t>c_{M+1}$
\[L([c_{M+1},t-1])\leq L([A,+\infty))\leq \epsilon.\]

By Lemma \ref{L:singvaluppereL}, 
\[\|R_{\alpha}\|_e=\inf_k a_k(R_\alpha) \leq a_{2N+2}(R_{\alpha})\leq \frac{\epsilon}{\sqrt{2}}.\]
In particular, $\|R_{\alpha}\|_e \leq \frac{\mathcal{L}}{\sqrt{2}}$.

To prove the reverse inequality, we may assume that $\mathcal{L}>0$, and let $\eta \in(0, \mathcal{L})$. Then, $L([a,b]) >\eta$ for all sufficiently large $a$ and $b$, which shows that the $(\eta,L)$-sequence is
infinite. Hence, by Lemma \ref{L:lowerboundaN}, $a_{N}(R_{\alpha})\geq \frac{\eta}{\sqrt{2}}$ for all $N\in\mathbb{N}$, and so $\|R_{\alpha}\|_e\geq\frac{\mathcal{L}}{\sqrt{2}}$.
\end{proof}

In analogy with Theorem \ref{T:boundRaeLseq}, we can state the following result.

\begin{theorem} \label{T:compRaeLseq2}
    The following conditions are equivalent.
    \begin{enumerate}[label=(\roman*)]
        \item $R_\alpha$ is a compact operator on $\ell^2$.
        \item $\mathcal{L}:=\lim_n L([n,+\infty))=0$.
        \item Every $(\epsilon,L)$-sequence has finite length.
    \end{enumerate}
\end{theorem}
\begin{proof}
    If $(i)$ holds, then by Theorem \ref{T:compRaeLseq1}, $(ii)$ holds as well. Also, it was shown in the proof of the same theorem that for every $\epsilon>\mathcal{L}$ the associated $(\epsilon,L)$-sequence has finite length. This fact proves that $(ii)\implies(iii)$. To conclude the proof, we show that $(iii)\implies(i)$. Let $\epsilon>0$. Then, if $N$ is the length of the associated $(\epsilon,L)$-sequence, by Lemma \ref{L:singvaluppereL},
    \[a_{2N+2}\leq \frac{\epsilon}{\sqrt{2}}.\]
    In particular, $\|R_\alpha\|_e\leq\frac{\epsilon}{\sqrt{2}}$ for every $\epsilon>0$, proving that $R_\alpha$ is compact.
\end{proof}

\section{Schatten class of $R_\alpha$} \label{S:Schatten}
So far, we have related the boundedness and the compactness of $R_\alpha$ to both the values $(\sigma_k)_k$  introduced in Definition \ref{D:defsigma} and to the $(\epsilon,L)$-sequences. Now, we use the $(\epsilon,L)$-sequences to provide information about the Schatten class properties of $R_\alpha$. In particular, we relate the decay of the approximation numbers of $R_\alpha$ to the decay of $(\sigma_k)_k$. We will need the next technical lemmas to relate the quantities $J(I)$ and $L(I)$ to the sequence $(\sigma_k)_k$. We start with the quantities $A_c,B_c$ introduced in Definition \ref{D:JAcBc}. We recall the notation $Z_k=[2^k,2^{k+1}-1].$

\begin{lemma}\label{lemma 17}
   Given any interval $I\in\mathbb{N}$ with endpoints $a$ and $b$, let  \[
   A(I)=\sup_{m\in I} \,(b-m)^{\frac{1}{2}}\biggl(\sum_{j=a}^{m} |\alpha_j|^2\biggr)^{\frac{1}{2}}, \qquad  B(I)=\sup_{m\in I} \,(m-a)^{\frac{1}{2}}\biggl(\sum_{j=m}^{b} |\alpha_j|^2\biggr)^{\frac{1}{2}}.\]
   Then, for all $k\in\mathbb{N}$,
\[
A(Z_{k}\cup Z_{k+1})\geq \sigma_{k}, \qquad B(Z_{k}\cup Z_{k+1})\geq \frac{\sigma_{k+1}}{2}.
\]
\end{lemma}
\begin{proof}
We notice that $Z_k \cup Z_{k+1}=[2^k, 2^{k+2}-1]$. On the one hand,
\begin{align*}
 A(Z_{k}\cup Z_{k+1})&= \sup_{2^k\leq m< 2^{k+2}-1}(2^{k+2}-1-m)^{\frac{1}{2}}\biggl(\sum_{j=2^k}^{m} |\alpha_j|^2\biggr)^{\frac{1}{2}}\\
 &\geq  (2^{k+2}-2^{k+1})^{\frac{1}{2}}\biggl(\sum_{j=2^k}^{2^{k+1}-1} |\alpha_j|^2\biggr)^{\frac{1}{2}}\\
 &\geq (2^{k+1})^{\frac{1}{2}}2^{-\frac{k+1}{2}}\sigma_{k}=\sigma_k.
\end{align*}
We have used \eqref{E:sigmabounds}. Similarly, for $B$ we have that
\begin{align*}
 B(Z_{k}\cup Z_{k+1})&= \sup_{2^k< m\leq 2^{k+2}-1}(m-2^{k})^{\frac{1}{2}}\biggl(\sum_{j=m}^{2^{k+2}-1} |\alpha_j|^2\biggr)^{\frac{1}{2}}\\
 &\geq  (2^{k+1}-2^{k})^{\frac{1}{2}}\biggl(\sum_{j=2^{k+1}}^{2^{k+2}-1} |\alpha_j|^2\biggr)^{\frac{1}{2}}\\
 &\geq (2^{k})^{\frac{1}{2}}2^{-\frac{k+2}{2}}\sigma_{k+1}=\frac{\sigma_{k+1}}{2}.
\end{align*}
\end{proof}

Using the notation of Definition \ref{D:JAcBc}, $A(I)=A_b^I$ and $B(I)=B_a^I$. 

The following three lemmas are essential to prove Theorem \ref{T:main3}. The first one shows that, if an interval $I$ is big enough so that there exist enough numbers $k\in\mathbb{N}$ such that $Z_k\subseteq I$ with $\sigma_k >\epsilon,$ then we have a lower bound for $J(I)$ that depends on $\epsilon$. The second lemma shows that the length of the $(\epsilon,L)$-sequence controls the amount of “large" $\sigma_k$'s. Finally, the third lemma shows that the number of large approximation numbers controls the number of large $\sigma_k$'s. 

\begin{lemma}\label{L:lemmaS(e)}
Let $I$ be an interval in $\mathbb{N}$ with endpoints $a$ and $b$. We consider
$J(I)$ as in Definition \ref{D:JAcBc}. Let $\epsilon > 0$ and suppose that  
\begin{equation}\label{E:S(e)}
S(\epsilon):=\{k\in \mathbb{N}: Z_k\subseteq I, \sigma_{k}>\epsilon\}
\end{equation}
has at least $4$ distinct elements. Then, $J( I ) > \dfrac{\epsilon}{2}$ and $L(I)>\dfrac{\epsilon}{2\sqrt{2}}$.
\end{lemma}
\begin{proof}
Let $k_1:=\min\{k\in\mathbb{N}\colon k\in S(\epsilon)\}$ and $k_2<k_3<k_4$ be elements in $S(\epsilon)\setminus\{k_1\}$, and let $c \in (a, b)$.  Notice that, either $Z_{k_1}\cup Z_{k_2} \subseteq [a,c)$, or $Z_{k_3}\cup Z_{k_4}\subseteq (c,b]$. Without loss of generality, we assume that the former holds. Then, $Z_{k_1} \cup Z_{k_1+1}\subseteq [a,c)$ and, by Lemma \ref{lemma 17},
\begin{align*}
    A_c&=A([a,c])=\sup_{m\in [a,c]}(c-m)^{\frac{1}{2}}\biggl(\sum_{j=a}^{m} |\alpha_j|^{2}\biggr)^{\frac{1}{2}}\\ 
    &\geq \sup_{m\in [a,2^{k_1+2}-1]}(2^{k_1+2}-1-m)^{\frac{1}{2}} \biggl(\sum_{j=2^{k_1}}^{m} |\alpha_j|^{2}\biggr)^{\frac{1}{2}}\\ 
    & \geq A(Z_{k_{1}}\cup Z_{k_1+1})\geq \sigma_{k_1} > \epsilon.
\end{align*}

A similar argument shows that, if $Z_{k_3}\cup Z_{k_4}\subseteq (c,b]$, then $Z_{k_4-1} \cup Z_{k_4}\subseteq (c,b]$ and, by Lemma \ref{lemma 17},
\begin{align*}
   B_c
    &\geq B(Z_{k_{4}-1}\cup Z_{k_4})\geq \frac{\sigma_{k_4}}{2}> \frac{\epsilon}{2}.
\end{align*}
Therefore,
\begin{align*}
 J(I)&=\inf_{c \in I}\max(A_c,B_c) > \frac{\epsilon}{2},  
\end{align*}
and, because of Lemma \ref{L:relationKJL},
\[
L(I)=\sqrt{2}K(I)\geq \frac{J(I)}{\sqrt{2}}>\frac{\epsilon}{2\sqrt{2}}.
\]
\end{proof}

We recall that $\#E$ denotes the cardinality of the set $E\subseteq \mathbb{N}$. 

\begin{lemma}\label{lemma 20}
    Let $t>0$, $\epsilon_t=\dfrac{t}{2\sqrt{2}}$ and let $N = N(\epsilon_t)$ be the length of the $(\epsilon_t,L)$-sequence $\{c_i\}_{i=0}^N$. Then,
\[\#\{k\in\mathbb{N}: \sigma_{k}>t\}\leq 5N(\epsilon_t)+3.\]
\end{lemma}
\begin{proof}
By definition of $(\epsilon_t,L)$-sequence, for $i=0,\ldots,N$ we have that  $L([c_i,c_{i+1}-2])\leq \epsilon_t.$
We introduce the set 
\begin{align*}
    E:=&\{k\in\mathbb{N}:\exists \,i\in \{1,\ldots,N\} \,\, \text{such that}\,\, c_i \in Z_{k} \,\,\text{or}\,\, c_i -1 \in Z_k\}.
\end{align*}
Notice that, for every $k\in\mathbb{N}$, $c_0=0\notin Z_k$. Since the intervals $Z_k$ are pairwise disjoint, $\# E \leq 2N$.

 On the other hand, if $k\in \mathbb{N}\setminus E$, then there exists $i\in \{0,\ldots,N\}$ such that $Z_k\subseteq [c_i+1,c_{i+1}-2].$ Then, we have that
 \[  \#\{k\in\mathbb{N}\setminus E: \sigma_{k}>t\} \leq 3(N+1).\]
This is true because for every $i\in\{0,\ldots,N\}$ there can be at most three different elements in $\{k\in\mathbb{N}\setminus E: \sigma_{k}>t\}$. In fact, if $i\in\{0,\ldots,N-1\}$ and we had four different $k_1,k_2,k_3,k_4$ such that $Z_{k_j}\subseteq [c_i,c_{i+1}-2]$ and $\sigma_{k_j}>t$ for $j=1,\ldots,4$, then by Lemma \ref{L:lemmaS(e)}, we would have that
\[
\epsilon_t=\frac{t}{2\sqrt{2}} < L([c_i,c_{i+1}-2])\leq \epsilon_t.  
\]
Also, if $i=N$, since every $Z_k$ has finite length, we would have $Z_{k_j}\subseteq [c_N,M]$ and $\sigma_{k_j}>t$ for $j=1,\ldots,4$, where $M$ is big enough. Consequently, for Lemma \ref{L:lemmaS(e)}, we would have that
\[\epsilon_t=\frac{t}{2\sqrt{2}} < L([c_N,M])\leq L([c_N+\infty))\leq \epsilon_t.\]
It follows that
\[\#\{k\in\mathbb{N}: \sigma_{k}>t\}\leq \#E+\#\{k\in\mathbb{N}\setminus E: \sigma_{k}>t\} \leq 5N+3.\]
\end{proof}

\begin{lemma}\label{lema 21}
Let $t>0$ be such that $\mathcal{L}<\epsilon_t=\dfrac{t}{2\sqrt{2}}\leq L(0,\infty)$.  Then,
\[\#\{k\in\mathbb{N}:\sigma_{k}>t\}\leq 5\#\left\lbrace k\in\mathbb{N}:a_k(R_{\alpha})> \frac{t}{8}\right\rbrace+3.\]
\end{lemma}
\proof 
Let $N(\epsilon_t)$ be the length of the $(\epsilon_t,L)$-sequence, that is finite as shown in the proof of Theorem \ref{T:compRaeLseq1}. By Lemma \ref{L:lowerboundaN}, for every $k\leq N(\epsilon_t)$,
\[a_k(R_\alpha)\geq a_{N(\epsilon_t)}(R_\alpha) \geq \frac{\epsilon_t}{\sqrt{2}}>\frac{\epsilon_t}{2\sqrt{2}}.\]
In particular,
\[\#\left\lbrace k\in\mathbb{N}:a_k(R_{\alpha})> \frac{\epsilon_t}{2\sqrt{2}}\right\rbrace \geq N(\epsilon_t).\]
By Lemma \ref{lemma 20}, we conclude that  
\[ \#\{k\in\mathbb{N}:\sigma_{k}>t\}
 \leq 5N(\epsilon_t)+3\leq 5\#\left\lbrace k\in\mathbb{N}:a_k(R_{\alpha})> \frac{\epsilon_t}{2\sqrt{2}}\right\rbrace+3. \qedhere\]

We are finally ready to explicitly relate the $\ell^p$-norm of the sequence $\sigma=(\sigma_k)_k$ to the $\ell^p$-norm of the approximation numbers $(a_n(R_\alpha))_n$. 
\begin{theorem}\label{lemma 22}
For any $p > 0$, there exists a constant $C_p$ depending only on $p$ such that
    \[\|(\sigma_{k})_k\|_{p}\leq C_p \|(a_{k}(R_\alpha))_k\|_{p}.\]
 Here, the $\ell^p$-(quasi)norms have their natural meanings.
\end{theorem}
\begin{proof}
We assume that $(a_{k}(R_\alpha))_k \in \ell^p$, otherwise the statement is trivial. In particular, $\lim_k a_k(R_\alpha)=0$, and therefore $R_\alpha$ is compact. Let $S = \sup_k|\sigma_{k}|$. Since $R_\alpha$ is bounded, by Theorem \ref{T:boundRasigma} we have that $S<\infty$. Then, from  Proposition 1.1.4 \cite{grafakos2014classical} and Lemma \ref{lema 21}, we have that
 \begin{align*}
\|\sigma\|^{p}_{p}&=p\int_{0}^{S}t^{p-1}\#\{k\in\mathbb{N}:\sigma_{k}>t\}\,\operatorname{d}\!t\\
&\leq 5p\int_{0}^{\infty}t^{p-1}\#\left\lbrace k\in\mathbb{N}:a_{k}(R_\alpha)> \frac{t}{8}\right\rbrace\,\operatorname{d}\!t+3S^p\\
&= 5p\, 8^p\int_{0}^{\infty}w^{p-1}\# \left\lbrace k\in\mathbb{N}:a_{k}(R_\alpha)> w\right\rbrace \,\operatorname{d}\!w+3S^p\\
&= 8^p5 \|(a_{k}(R_\alpha))_k\|_{p}^{p}+3S^{p}.
 \end{align*} 
Since, by Theorem \ref{T:boundRasigma}, $S=\|\sigma\|_{\infty}\leq \sqrt{2}\|R_\alpha\|$, we have that
\begin{align*}
\|\sigma\|_{p}\leq& \Big( 8^p5 \|(a_{k}(R_\alpha))_k\|_{p}^{p}+2^{\frac{p}{2}}3 a_1(R_\alpha)^p\Big)^{1/p}\\
\leq& \|(a_{k}(R_\alpha))_k\|_{p} \big(  8^p5 +2^{\frac{p}{2}}3\big)^{1/p}.    
\end{align*}
\end{proof}
We obtain the reverse inequality of Theorem \ref{lemma 22}. First, a technical lemma. We recall that $Z_k=[2^k,2^{k+1}-1]$.

\begin{lemma}\label{L:k0k1k2}
Let $j_0, j_2 \in \mathbb{N}$ with $j_0 +1< j_2$, and let $I = [a, b]$ be an interval in $\mathbb{N}$ such that $a\geq 2^{j_0}, b\in Z_{j_{2}}$. Then, for every $s\in I$,
\begin{align}
S:= (s-a)^{\frac{1}{2}}\biggl(\sum_{i=s}^b|\alpha_i|^2\biggr)^{\frac{1}{2}}\leq 2\max_{j_{0}\leq j\leq j_{2}}\sigma_{j}. 
\end{align}
\end{lemma}

\begin{proof}
We have to distinguish three cases: when $s\in Z_{j_0}$, when $s\in Z_{j_2}$, and when $s\in I\setminus (Z_{j_0}\cup Z_{j_2}).$ In the first case, $a\leq s\leq 2^{j_{0}+1}-1$, and by \eqref{E:sigmabounds} we have that
\begin{align*}
(s-a)\sum_{i=s}^{b} |\alpha_i|^2&\leq 2^{j_{0}+1} \sum_{j=j_{0}}^{j_{2}} \sum_{i=2^{j}}^{2^{j+1}-1}|\alpha_i|^2\\
&\leq 2^{j_{0}+1} \sum_{j=j_{0}}^{j_{2}} 2^{-j} \sigma_j^2 \leq 4\max_{j_{0}\leq j \leq j_{2}} \sigma_j^2.
\end{align*}
In the second case, that is, $2^{j_{2}}\leq s\leq b$, similarly,
\begin{align*}
(s-a)\sum_{i=s}^{b} |\alpha_i|^2&\leq 2^{j_{2}+1} \sum_{i=2^{j_{2}}}^{2^{j_{2}+1}-1}|\alpha_i|^2 \leq 2 \sigma_{j_{2}}^2.
\end{align*}
Finally, in the last case, there exists $j_1\in\mathbb{N}$ such that $j_0<j_1<j_2$ and $s\in Z_{j_1}$. Then, 
\begin{align*}
S^2&\leq (2^{j_{1}+1}-1-2^{j_{0}}) \sum_{k=2^{j_{1}}}^{2^{j_2+1}-1}|\alpha_k|^2 \leq 2^{j_{1}+1}\sum_{n=j_{1}}^{j_{2}}\frac{\sigma_{n}^{2}}{2^{n}} \leq 4\max_{j_1\leq j\leq j_2} \sigma_j^2. 
\end{align*}
The result now follows.
\end{proof}
We recall that
\[\mu(I)=\sum_{k\in I} |\alpha_k|^2.\]
\begin{theorem}\label{Estim 1}
   For any $p\in (1,+\infty)$, there exists a positive constant $C_p$ depending only on $p$ such that
    \[
    \|(a_{k}(R_\alpha))_k\|_{p} \leq C_p \left(\sum_{k \geq -1}|\sigma_k|^p\right)^{1/p}.
    \]
\end{theorem}

\begin{proof}
Let $\epsilon>0$ and $I_k = [c_{k}, c_{k+1}-1]$ be the intervals associated to the $(\epsilon,L)$-sequence $c_k$. We introduce the families
\[\mathcal{A}_j=\{I_k\colon I_k\subseteq \{0\}\cup Z_j\cup Z_{j+1}\},\]

for any $j \in \mathbb{N}$. We notice that $\mathcal{A}_0$ contains all the $I_k$'s in $[0,3]$, if any. In particular, $\#\mathcal{A}_0\leq 2$. According to Lemma \ref{L:relationKJL} and Definition \ref{D:JAcBc}, we notice that
\begin{align*}
 \epsilon \#\mathcal{A}_j&\leq 2\sqrt{2} \sum_{k\colon I_k\in\mathcal{A}_j} J([c_k,c_{k+1}-1])\\
 &\leq 2\sqrt{2}\sum_{k\colon I_k\in\mathcal{A}_j} \max\{A_{c_k}^{I_k},B_{c_k}^{I_k}\} \\
 &=2\sqrt{2}\sum_{k\colon I_k\in\mathcal{A}_j} \sup_{c_k\leq s\leq c_{k+1}-1} (s-c_k)^\frac{1}{2}\biggl(\sum_{i=s}^{c_{k+1}-1} |\alpha_i|^2 \biggr)^\frac{1}{2}\\
 &\leq 2\sqrt{2} \sum_{k\colon I_k\in\mathcal{A}_j} (c_{k+1}-1-c_k)^\frac{1}{2}\mu(I_k)^\frac{1}{2}.
 \end{align*}
 Consequently, by Cauchy-Schwarz and \eqref{E:sigmabounds}, for $j \geq 1$, we have that
 \begin{align*}
 \epsilon \#\mathcal{A}_j &\leq  2\sqrt{2} \biggl(\sum_{k\colon I_k\in\mathcal{A}_j}(c_{k+1}-1-c_k)\biggr)^\frac{1}{2}   \biggl(\sum_{k\colon I_k\in\mathcal{A}_j}\mu(I_k)\biggr)^\frac{1}{2} \\
 &\leq 2\sqrt{2} \biggl(2^{j+2}-1-2^j\biggr)^\frac{1}{2}\mu(Z_j\cup Z_{j+1})^\frac{1}{2}\\
  &\leq 2\sqrt{2} \biggl(2^{j}3 
  \sum_{i=2^j}^{2^{j+2}-1}|\alpha_i|^2\biggr)^\frac{1}{2} \\
&\leq 2\sqrt{6}(\sigma_j^2 + \frac{1}{2}\sigma_{j+1}^2)^\frac{1}{2}.
\end{align*}
Thus, for $j\geq 1$ we have that \[\#\mathcal{A}_j\leq n_j:= [6 \max\left(\sigma_{j}, \sigma_{j+1}\right)/ \epsilon],\] where $[x]:=\max\{n\in\mathbb{Z}\colon n\leq x \}$ is the integer part of the real number $x$. Therefore, we have that
\begin{align*}
    \sum_{j\in\mathbb{N}}\#\mathcal{A}_j &\leq 2+\sum_{j\geq 1}\sum_{i=1}^{n_j} 1 \\
    &=2+ \sum_{i=1}^{\max_j n_j}\#\{j\in\mathbb{N}\colon n_j \geq i\} \\
    &\leq2+ \sum_{i=1}^\infty \#\{j\in\mathbb{N}\colon6\max\left( \sigma_{j}, \sigma_{j+1}\right)/ \epsilon \geq i\} \\
    &\leq 2+2\sum_{i=1}^\infty \#\{k\in\mathbb{N}\colon  6 \sigma_{k} /\epsilon \geq i\},
\end{align*}
where the constant $2$ in the last inequality derives from the fact that two consecutive $j$'s may have the same $\max(\sigma_j,\sigma_{j+1})$.
On the other hand, let $\{I_k\}_{k\in \mathcal{Q}}=\{I_k\colon I_k \notin \bigcup_{j\in\mathbb{N}}\mathcal{A}_j\}$. For each $k\in\mathcal{Q}$ there exist $j_{0,k}<j_{1,k}<j_{2,k}$ such that $c_k\in Z_{j_{0,k}}$ and $c_{k+1}-1\in Z_{j_{2,k}}$. Then, arguing as above, for every such $k\in\mathcal{Q}$,
\begin{align*}
    \epsilon&\leq  2\sqrt{2} \sup_{c_k\leq s\leq c_{k+1}-1} (s-c_k)^\frac{1}{2}\biggl(\sum_{i=s}^{c_{k+1}-1} |\alpha_i|^2 \biggr)^\frac{1}{2}.
\end{align*}
Let $s_k\in I_k$ be such that
\[ (s_k-c_k)^\frac{1}{2}\biggl(\sum_{i=s_k}^{c_{k+1}-1} |\alpha_i|^2 \biggr)^\frac{1}{2}=\sup_{c_k\leq s\leq c_{k+1}-1} (s-c_k)^\frac{1}{2}\biggl(\sum_{i=s}^{c_{k+1}-1} |\alpha_i|^2 \biggr)^\frac{1}{2}.\]

Then, by Lemma \ref{L:k0k1k2}, 
\[(s_k-c_k)^\frac{1}{2}\biggl(\sum_{i=s_k}^{c_{k+1}-1} |\alpha_i|^2 \biggr)^\frac{1}{2} \leq 2 \max_{j_0\leq n\leq j_2}\sigma_n.\]
We have shown that, for every $k\in\mathcal{Q}$, there exists $j_{0,k}\leq n_k\leq j_{2,k}$ such that
\[1\leq \frac{4\sqrt{2} \sigma_{n_k}}{\epsilon}.\]
 Also, notice that the same $\sigma_{n_k}$ my be obtained by at most two consecutive $k\in\mathcal{Q}$, showing that the mapping $k\mapsto n_k$ is “quantitatively almost injective".
Thus,
\[\#\mathcal{Q}=\#\{k\in\mathcal{Q}\colon 1\leq \frac{4\sqrt{2} \sigma_{n_k}}{\epsilon}\}\leq 2\#\{n\geq1\colon \frac{4\sqrt{2} \sigma_{n}}{\epsilon}\geq 1\}.\]
We obtain the following upper bound for the length of the $(\epsilon,L)$-sequence:
\begin{align} 
    \nonumber N(\epsilon)&\leq 2+\sum_{j\geq 1}\#\mathcal{A}_j+\#\mathcal{Q} \\
   \nonumber  &\leq 2+ 4 \sum_{n=1}^{\infty}\#\left\lbrace k\geq1:\sigma_{k}\geq \frac{n\epsilon}{6}\right\rbrace\\
   \label{E:Testim1}  &\leq  2+ 4 \sum_{n=1}^{\infty}\#\left\lbrace k\geq1:\sigma_{k}>\frac{n\epsilon}{6\sqrt{2}}\right\rbrace.
\end{align}
Since, by Lemma \ref{L:singvaluppereL},
\[
\{k\in\mathbb{N}:a_{k}(R_\alpha)>\frac{\epsilon}{\sqrt{2}}\}\subseteq[0, 2N(\epsilon)+2],
\]
we obtain that
\begin{align*}
 \|(a_{k}(R_\alpha))_k\|_{p}^{p}&=p\int_{0}^{\infty}t^{p-1}\#\{k\geq1\colon a_{k}(R_{\alpha})>t\}\,\operatorname{d}\!t  \\
 &\leq p\int_{0}^{\|R_{\alpha}\|}t^{p-1}\big(2N(\sqrt{2}t)+3\big)\,\operatorname{d}\!t.
 \end{align*}
 Now, by \eqref{E:Testim1} and an appropriate change of variable,
 \begin{align*}
     \int_{0}^{\|R_{\alpha}\|}\!\!\!t^{p-1}&N(\sqrt{2}t)\operatorname{d}\!t \leq \int_{0}^{\|R_{\alpha}\|}\!\!\!t^{p-1}\!\bigg(\!2+4\!\sum_{n=1}^\infty\!\#\biggl\lbrace k\in\mathbb{N}:\sigma_{k}>\frac{nt\sqrt{2}}{6\sqrt{2}} \biggr\rbrace\!\bigg)\!\operatorname{d}\!t \\
     &\leq \frac{2\|R_\alpha\|^p} {p}+ 
     6^{p}4 \sum_{n=1}^\infty \frac{1}{n^p}\int_0^\infty s^{p-1} \#\lbrace k\in\mathbb{N}:\sigma_{k}> s \rbrace\operatorname{d}\!s \\
     &\leq \frac{2\|R_\alpha\|^p}{p} +
     \frac{6^{p}4}{p} \bigg(\sum_{n=1}^\infty \frac{1}{n^p}\bigg)\|\sigma\|_p^p.
 \end{align*}
 We may assume that $(\sigma_k)_k\in\ell^p$, otherwise the inequality that we are trying to prove would be trivial. In particular, by Theorem \ref{T:boundRasigma},
 \[\|R_{\alpha}\| \leq 4\sqrt{2}\sup_{k\geq -1} |\sigma_k| \leq 4\sqrt{2}\left(\sum_{k \geq -1}|\sigma_k|^p\right)^{1/p}\]
 and then we conclude that
 \begin{align*}
 \|(a_{k}(R_\alpha))_k\|_{p}^{p}&\leq 3\|R_\alpha\|^p+4 \|R_\alpha\|^p+\left(6^{p}8 \sum_{n \geq 1}\frac{1}{n^p}\right)\|\sigma\|_p^p  \leq C_p^p \|\sigma\|_p^p.
 \end{align*}
\end{proof}

As a consequence of Theorems \ref{lemma 22} and \ref{Estim 1}, we have the following characterization for when $R_\alpha$ belongs in a Schatten class. 

\begin{theorem}\label{T:Schattensigma}
    Let $\alpha\in\ell^2$ and $1<q<\infty$. The operator $R_\alpha$ belongs to the Schatten class $\mathcal{S}^q$ if and only if the sequence $(\sigma_k)_k$ belongs to $\ell^q$.
\end{theorem}

We are ready to prove the third main theorem. 

\begin{proof}[Proof of Theorem \ref{T:main3}]
We recall that $f_\alpha(z)=\sum_n \alpha_n z^n$ and that
\[\Delta_n f_\alpha(z)=\sum_{j=2^n}^{2^{n+1}-1}\alpha_j z^j.\]
The equivalence of \ref{T:RaSchattenii} and \ref{T:RaSchatteniii} was already discussed in Section \ref{S:mainresults}. 

By \eqref{E:sigmabounds}, condition \ref{T:RaSchatteniii} holds if and only if $\sigma\in\ell^q$, for
\begin{align*}
   \sum_{n}\bigl( 2^n \|\Delta_nf_\alpha\|_{H^2}^2\bigr)^{\frac{q}{2}}&=\sum_{n}\biggl( 2^n \sum_{j=2^n}^{2^{n+1}-1}|a_j|^2\biggr)^{\frac{q}{2}}\leq \|\sigma\|_{\ell^q}^q
\end{align*}
and
\begin{align*}
   \|\sigma\|_{\ell^q}^q&\leq \sum_{n}\biggl( 2^{n+1} \sum_{j=2^n}^{2^{n+1}-1}|a_j|^2\biggr)^{\frac{q}{2}}\leq 2^{\frac{q}{2}}\sum_{n}\bigl( 2^n \|\Delta_nf_\alpha\|_{H^2}^2\bigr)^{\frac{q}{2}}.
\end{align*}
The proof now follows from Theorem \ref{T:Schattensigma}.
\end{proof}

 \section{Proof of Theorem \ref{T:main4}} \label{S:section6}
In the previous sections we have considered the elements of $\ell^2$ as complex-valued functions defined on $\mathbb{N}$. For what follows, we will instead work with infinitely long vectors $x=(x_0,x_1,\ldots)$.
We recall that the shift operator acting on $\ell^2$ is defined as
\[S \colon \ell^2 \ni (x_0,x_1,x_2,\ldots)\mapsto(0,x_0,x_1,\ldots)\in\ell^2\]
and its adjoint, the backward shift, is the operator
\[S^*\colon\ell^2\ni (x_0,x_1,x_2,\ldots)\mapsto(x_1,x_2,x_3,\ldots)\in \ell^2.\]
They satisfy the relation $SS^*=I-\textbf{e}^0 \otimes \textbf{e}^0$, where $\textbf{e}^0 \otimes \textbf{e}^0$ is the rank-one operator that assigns to each $x\in\ell^2$ the sequence $(x_0,0,0,\ldots)$.

We begin with a lemma concerning the eigenvalues of the matrix $T_c$. This result already appeared in \cite{DellepianeCesaro}.

\begin{lemma} \label{L:eigenvaluesTc}
Let $(c_{k})_{k\in\mathbb{N}}$ be a sequence of complex numbers. Then $c_k$, $k \geq 0$, is an eigenvalue of the matrix $T_c$, with corresponding eigenvector
\[\mathbf{v}_k =\sum_{i=0}^{k} \mathbf{e}^i.\]
\end{lemma}

\begin{proof}
We prove by induction that $T_c \mathbf{v}_k = c_k \mathbf{v}_k$. For $k=0$, the result is trivial. Assume that it holds for $k\in\mathbb{N}$. Then
\begin{align*}
T_c \mathbf{v}_{k+1} &= T_c \mathbf{v}_k + T_c \mathbf{e}_{k+1} = c_k \mathbf{v}_k + \sum_{i=0}^k (c_{k+1}-c_{k})\mathbf{e}_{i} + c_{k+1} \mathbf{e}_{k+1} \\
&= c_k \mathbf{v}_k + (c_{k+1}-c_{k}) \mathbf{v}_k +c_{k+1}\mathbf{e}_{k+1} = c_{k+1}\mathbf{v}_{k+1}. \qedhere
\end{align*}
\end{proof}

Let $(c_k)_k$ be a sequence of complex numbers. We introduce the Rhaly matrix $A_c$ defined in the following way:
    \[A_c:=
\begin{pmatrix}
c_1-c_0 & 0 & 0 & 0 & 0 & \cdots\\
c_2-c_1 & c_2-c_1 & 0 & 0 & 0 & \cdots\\
c_3-c_2 & c_3-c_2 & c_3-c_2 & 0 & 0 & \cdots \\
c_4-c_3 & c_4-c_3 & c_4-c_3 & c_4-c_3 & 0 & \cdots\\
\vdots & \vdots & \vdots & \vdots &\vdots & \ddots
\end{pmatrix}.
\]

\begin{lemma} \label{L:Ac1}
The following properties hold.
\begin{enumerate}
    \item $T_c$ is bounded on $\ell^2$ if and only if $(c_k)_k$ is bounded and $A_c$ is bounded on $\ell^2$.
    \item $T_c$ is compact on $\ell^2$ if and only if
   $\lim_n c_k=0$ and $A_c$ is compact.
    \item For $1<q<\infty$, $T_c$ belongs to $\mathcal{S}^q(\ell^2)$ if and only if $(c_k)_k\in\ell^q$ and $A_c\in\mathcal{S}^q(\ell^2)$.
\end{enumerate}
\end{lemma}
\begin{proof}
We will work with the adjoint $(T_c)^*$. Without loss of generality, we can omit the conjugation. We can write $(T_c)^*$ as the sum of two operators $L_c + D_c$, where
\[
L_c=\begin{pmatrix}
 0 & 0 & 0 &  \cdots\\
c_1-c_0 & 0 & 0 &  \cdots\\
c_2-c_1 & c_2-c_1 & 0 &   \cdots\\
c_3-c_2 & c_3-c_2 & c_3-c_2 &   \cdots \\
\vdots & \vdots & \vdots &  \ddots
\end{pmatrix},\,\, D_c= \begin{pmatrix}
c_0 & 0 & 0 & 0 & \cdots\\
0 & c_1 & 0 & 0 & \cdots\\
0 & 0 & c_2 & 0 & \cdots\\
0 & 0 & 0 & c_3 & \cdots \\
\vdots & \vdots & \vdots & \vdots & \ddots
\end{pmatrix} .
\]
Notice that the composition with the backward shift operator corresponds to erasing the first row of the matrix:
\[S^* L_c = \begin{pmatrix}
c_1-c_0 & 0 & 0 & 0 & 0 & \cdots\\
c_2-c_1 & c_2-c_1 & 0 & 0 & 0 & \cdots\\
c_3-c_2 & c_3-c_2 & c_3-c_2 & 0 & 0 & \cdots \\
c_4-c_3 & c_4-c_3 & c_4-c_3 & c_4-c_3 & 0 & \cdots\\
\vdots & \vdots & \vdots & \vdots &\vdots & \ddots
\end{pmatrix} = A_c.\]
In light of the relation $SS^*=I-\textbf{e}^0 \otimes \textbf{e}^0$, we have that
\[S(S^* L_c)x = L_c x - (L_c x)_0 \textbf{e}^0 = L_c x, \qquad x\in \ell^2.\]
Thus, $L_c = S A_c$. It is known that boundedness of the diagonal matrix $D_c$ is equivalent to the boundedness of $(c_k)_k$, that its compactness is equivalent to the limit of $(c_k)_k$ being $0$, and that $D_c\in\mathcal{S}^q(\ell^2)$ if and only if $(c_k)_k\in\ell^q$. This proves all the three sufficient conditions of Lemma \ref{L:Ac1}.
We move to the necessary conditions. If $T_c$ is bounded, then the sequence $(c_k)_k$ is bounded as well, since they are all eigenvalues of $T_c$ by Lemma \ref{L:eigenvaluesTc}. Then, 
\[A_c = S^*((T_c)^*-D_c)\]
is bounded. The compactness follows from the same argument: since every $c_k$ is an eigenvalue of $T_c$, if $T_c$ is compact, then the sequence $(c_k)_k$ can accumulate only at $0$.  Finally, if $T_c\in\mathcal{S}^q(\ell^2)$, then 
\[\sum_k |\langle T_c \textbf{e}^k,\textbf{e}^k\rangle|^q = \sum_k |c_k|^q<\infty,\]
showing that $D_c\in\mathcal{S}^q(\ell^2)$ and then $A_c\in\mathcal{S}^q(\ell^2)$ as well.
\end{proof}
We conclude proving the fourth main theorem.

\begin{proof}[Proof of Theorem \ref{T:main4}]
    The proof follows from Lemma \ref{L:Ac1} and, respectively,  Theorems \ref{T:main1}, \ref{T:main2} and \ref{T:main3}, applied to the Rhaly matrix $A_c$.
\end{proof}

\section{Conclusion}
The $(\epsilon,L)-$sequence technique may be useful for studying approximate values of other significant operators. First of all, one shall consider the factorable matrices, already defined in \eqref{factorable matrix}, which depend on two sequences $\alpha_n$ and $\beta_n$. The arguments presented in this article may adapt also to these operators. In this situation, the conditions on $f_\alpha$ will be in terms of weighted Besov spaces. 

There are still several open questions related to the Rhaly matrices. First of all, Theorem \ref{lemma 22} only provides a necessary condition for the $q$-Schatten class membership, for $0<q\leq1$. A complete characterization is missing.

The study of Rhaly matrices acting between conformally invariant spaces is interesting and techniques similar to \cite{Bellavita_Stylogiannis_2024} may be helpful. 

The $L$-shaped matrices\footnote{Not to be confused with the $L$-matrices!} are also related to the Rhaly matrices. Given a sequence $(\ell_n)_n$, the matrix
\[
L_= \begin{pmatrix}
\ell_0  & \ell_1 & \ell_2 &  \cdots\\
\ell_1 & \ell_1 & \ell_2 & \cdots\\
\ell_2 & \ell_2 & \ell_2 & \cdots \\
\vdots & \vdots & \vdots & \ddots
\end{pmatrix} 
\]
is called $L$-shaped.
They are commonly used in the study of Hadamard multipliers in function spaces and, recently, their spectral properties have caught a lot of attention, see \cite{Bouthat2021}, \cite{BOUTHAT20221} and \cite{STAMPACH2022109401} as examples. It is well known that, if $R_\alpha$ is a Rhaly matrix, then 
\[L:=R^*_\alpha R_\alpha\]
is an $L$-shaped matrix. It is reasonable to guess that some properties of the approximate values of $L$ may be obtained by using $(\epsilon,L)-$sequences related to $\alpha$. We will pursue these ideas in future works.

\bibliographystyle{plain}
\bibliography{dirichletmatrix}

\end{document}